\documentclass[12pt]{amsart}
\usepackage{amsmath,amssymb,amsfonts,amscd,verbatim,graphicx}
\usepackage{xypic}
\newtheorem{theorem}[subsection]{Theorem}
\newtheorem{lemma}[subsection]{Lemma}
\newtheorem{proposition}[subsection]{Proposition}

\theoremstyle{definition}
\newtheorem{definition}[subsection]{Definition}
\newtheorem{example}[subsection]{Example}

\theoremstyle{remark}

\numberwithin{equation}{section}

\theoremstyle{definition}

\newcommand{\co}{\colon\thinspace}

\newcommand{\nl}{\hfil\break}

\begin{document}

\title[Link groups of $4$-manifolds]{Link groups of ${\mathbf 4}$-manifolds}
\author{Vyacheslav Krushkal}
\address{Department of Mathematics, University of Virginia, Charlottesville, VA 22904}
\email{krushkal\char 64 virginia.edu}

\thanks{This research was partially supported by the NSF}

\begin{abstract}
The notion of a {\em Bing cell} is introduced, and it is used to define
invariants, {\em link groups}, of $4$-manifolds. Bing cells
combine some features of both surfaces and $4$-dimensional handlebodies, and
the link group ${\lambda}(M)$ measures certain aspects of the handle
structure of a $4$-manifold $M$. This group is a quotient of the
fundamental group, and examples of manifolds are given with
${\pi}_1(M)\neq {\lambda}(M)$.
The main construction of the paper is a generalization of the Milnor group, which is
used to formulate an obstruction to embeddability of Bing
cells into $4$-space. Applications to the A-B slice problem and to the structure of topological arbiters are
discussed.
\end{abstract}

\maketitle

\section{Introduction} \label{introduction}
Maps of surfaces and of more general $2$-complexes have
been classically used to define invariants of topological spaces, for example
the fundamental group and the first homology group of a space. More
generally, one gets the quotients of ${\pi}_1 X$ by the terms of its
lower central series if one considers based loops
in a space $X$ modulo loops bounding maps of certain special
$2$-complexes, gropes \cite{FQ}. From this perspective gropes interpolate between
surfaces (null-homology) and disks (null-homotopy).

This paper introduces the notion of a {\em Bing cell}, which may be
viewed as a geometric dual to a grope.
The origin of this construction is in Milnor's theory of link
homotopy \cite{Milnor1}.
The idea in the definition of a Bing cell is to treat a collection of $2$-handles attached
to a homotopically essential link on an equal footing with an actual
$4$-dimensional $2$-handle $D^2\times D^2$,
see section \ref{outline sec} for a more detailed outline of the construction.
A Bing cell is not a $2$-complex, rather it is
a $4$-dimensional handlebody with a $2$-dimensional spine where the $4$-dimensional
thickening plays an important role.

In an analogy with the fundamental
group, the {\em link group} ${\lambda}(M)$ of a
$4$-manifold $M$ is defined as based loops in $M$ modulo loops bounding
Bing cells. The resulting invariant ${\lambda}(M)$ reflects
certain aspects of the handle structure of a $4$-manifold $M$, and it is not correlated with
the homology group $H_1(M)$.
Although their definition makes sense in any dimension, the link
groups are a non-trivial theory only in dimension $4$, and they are a topological but
{\em not} in general a homotopy invariant of a $4$-manifold.

Since the $2$-cell $D^2$ is a trivial example of a Bing cell, the link group
${\lambda}(M)$ is a quotient of the fundamental group of $M$. One
can easily find examples of $4$-manifolds with ${\pi}_1(M) \neq
{\lambda}(M)$. Bing cells may be given geometric and algebraic
gradings: {\em height} and {\em nilpotent class}, leading to a
two-parameter collection of link groups ${\lambda}_{i,j}(M)$, where
$j>i$. In this notation, the group ${\lambda}(M)$ above corresponds
to ${\lambda}_{\infty,\infty}$. We show that given a surjection of
finitely presented groups ${\pi}\longrightarrow {\lambda}$, where
${\pi}$ is aspherical, there are $4$-manifolds $M$ with
${\pi}_1(M)\cong{\pi}$, ${\lambda}_{1,2}(M)\cong{\lambda}$.

This work is motivated in part by the question of whether there
is a ``non-abelian'' Alexander duality in dimension $4$. This
question arises in the analysis of decompositions of the $4$-ball in
the $A,B$-slice problem, a reformulation of the $4$-dimensional
topological surgery conjecture. An application of link groups in
this context is given in \cite{K2}, and it is summarized in section \ref{applications sec} below.
For another recent application of the theory developed here, to the structure of topological arbiters,
see \cite{FK} and theorem \ref{uncountable} below.

The connection to the $A,B$-slice problem is provided by the
following theorem which is the main result of this paper, showing how Bing cells fit in the framework of
Milnor's theory of link homotopy:

\begin{theorem} \label{homotopically trivial} {\sl If the components of
a link $L\subset S^3=\partial D^4$ bound disjoint Bing cells in
$D^4$ then $L$ is homotopically trivial.}

\end{theorem}

This result is based on a generalization of the Milnor
group which is developed here in order to formulate an obstruction to embeddability of a disjoint
collection of Bing cells in $4$-space. We will next give a brief outline of the ideas underlying the
construction of Bing cells and of the generalized Milnor group, for more details see sections \ref{definitions section},
\ref{generalized Milnor}.

\subsection{Outline of the construction} \label{outline sec}

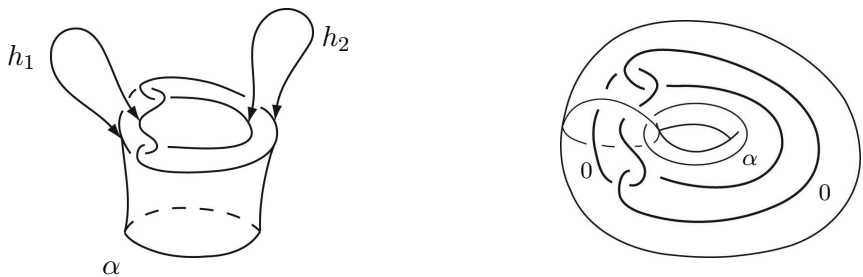
\begin{figure}[ht]
 \includegraphics[width=3.5cm]{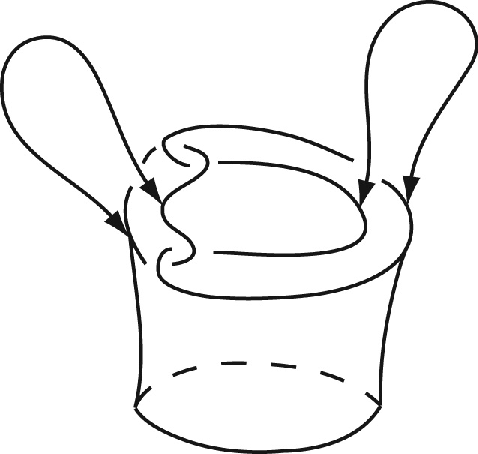} \hspace{3cm} \includegraphics[width=4cm]{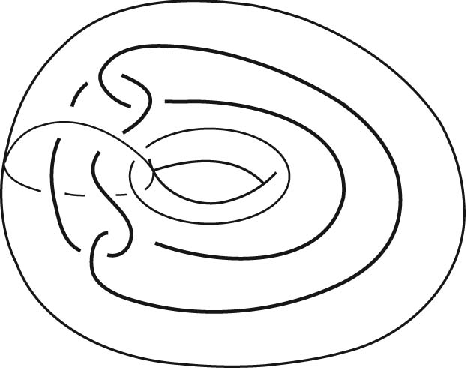}
{\small
    \put(-288,-6){${\alpha}$}
    \put(-324,73){$h_1$}
    \put(-205,80){$h_2$}}
{\scriptsize
    \put(-107,30){$0$}
    \put(-17,22){$0$}
    \put(-46,35){${\alpha}$}}
\caption{A pair of $2$-handles attached to the Bing double.}
\label{fig:surface Bing double}
\end{figure}

Consider the $4$-dimensional $2$-handle $H=D^2\times D^2$, thought of as a $4$-dimensional
thickening of its core $D^2\times\{ 0\}$. Remove a small disk $D^2_{\epsilon}$ from the core,
and consider the corresponding thickening $H_{\epsilon}=(D^2\smallsetminus D^2_{\epsilon})\times D^2$.
$H_{\epsilon}$ has a new part of the boundary, $S_{\epsilon}^1\times D^2$, which is the boundary
of $D^2_{\epsilon}\times D^2$ that was removed from the handle $H$. Attach to $H_{\epsilon}$ a pair of
zero-framed $2$-handles $h_1, h_2$ along the Bing double of the core of the solid torus $S_{\epsilon}^1\times D^2$, see figure \ref{fig:surface Bing double}.

This is the basic operation used in the construction of Bing cells, and it may be roughly
described as ``puncturing'' a $2$-handle and plugging in the puncture with $2$-handles whose
attaching curves form an {\em essential} link in the boundary of the puncture. Here the term ``essential''
is understood in the context of Milnor's theory of link homotopy, see \cite{Milnor1}
and section \ref{link-homotopy} in this paper. The most important example (motivating
the term {\em Bing cell})
is the Bing double, or more generally an iterated Bing double, of the core of the solid torus.

\begin{figure}[ht]
\centering
\bigskip
\includegraphics[width=5.75cm]{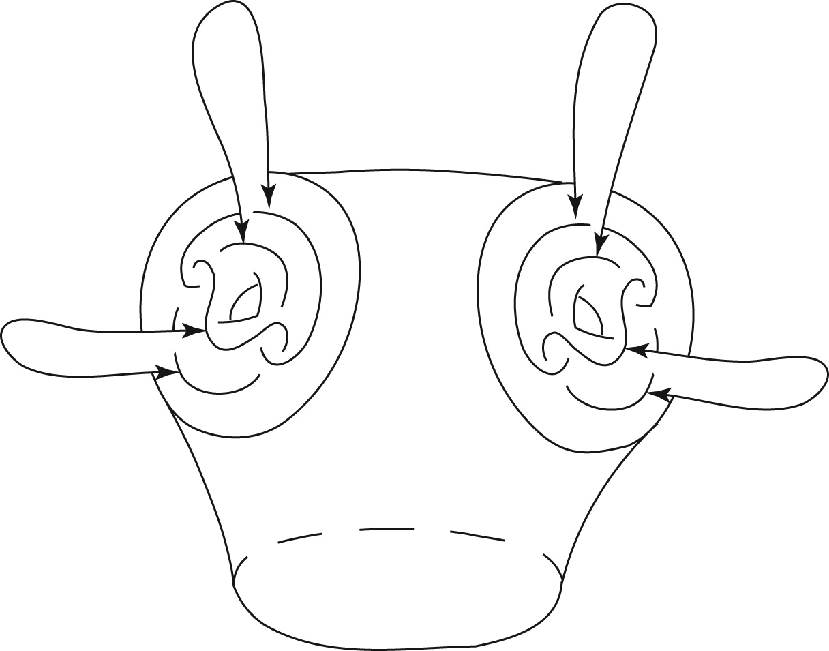}\hspace{2cm} \includegraphics[bb=0 -55 250 100, width=5.5cm]{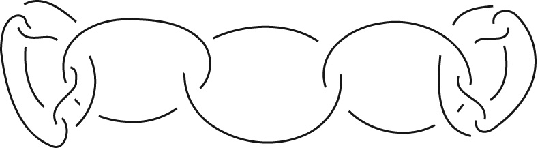}
{\small
    \put(-349,9){${\gamma}$}
    \put(-387,73){$h_1$}
    \put(-256,20){$C$}
    \put(-360,122){$h_2$}
    \put(-282,120){$h_3$}
    \put(-227,65){$h_4$}}
{\scriptsize
    \put(-164,58){$0$}
    \put(-141,80){$0$}
    \put(-18,79){$0$}
    \put(7,63){$0$}}
    {\small
    \put(-75,25){${\gamma}$}}
    \put(-118,73){\circle*{3}}
    \put(-39,72){\circle*{3}}
    \vspace{.3cm}
\caption{The handlebody $C$ (a Bing cell of height $1$): a schematic
picture and a Kirby diagram.}
\label{cell figure}
\end{figure}
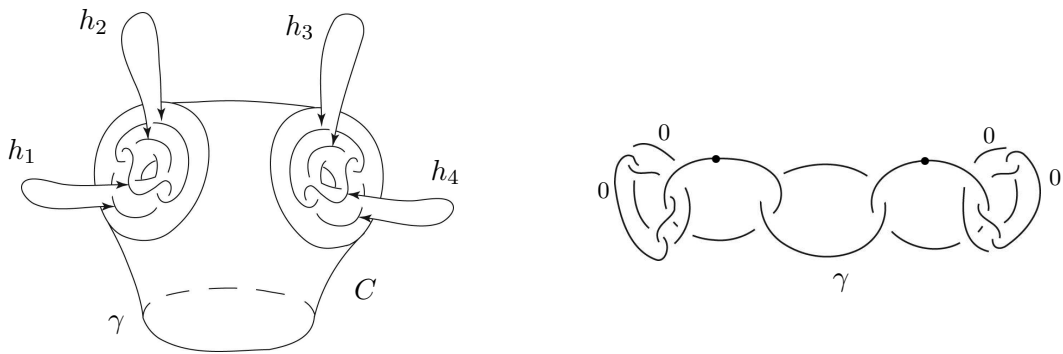

Now a {\em Bing cell of height $1$} is obtained from the $2$-handle $H$ by performing this operation
in a finite number of distinct locations in the core $D^2\times\{ 0\}$. For example, consider the case of two punctures
in more detail.
Let $P$ denote the pair of pants with boundary components ${\gamma}, {\alpha}_1, {\alpha}_2$, and let
$C$ denote $(P\times D^2)\; \cup$ four zero-framed $2$-handles $h_1\ldots,h_4$ attached to the Bing doubles of the curves ${\alpha}_1, {\alpha}_2$,
figure \ref{cell figure}. The operation above also may be applied to the handles $h_i$, leading to the construction of
a Bing cell of height $2$. Iterating this procedure (a finite number of times) yields an inductive construction of a general Bing cell.
There is a distinguished curve $\gamma$ in the boundary of a Bing cell $C$ which is the attaching curve of the original
$2$-handle $H$, ${\gamma}={\partial D^2}\times \{ 0\}$, and the term Bing cell will refer to a pair $(C,{\gamma})$.

The defining property in Milnor's theory of link homotopy \cite{Milnor1} is that the link components
are allowed to move by a homotopy so that {\em different components stay disjoint} from each other. It is natural to consider Bing cells
mapped into $4$-manifolds, requiring that the handles attached to different components of the Bing doubles are disjoint in the image. For example, in the case of the Bing cell $C$ in figure \ref{cell figure}, for a map $f\co C\longrightarrow M^4$ the requirement is that
$f(h_1)$ is disjoint from $f(h_2)$ and similarly $f(h_3)$ is disjoint from $f(h_4)$.

It is a classical fact \cite{Milnor1, Giffen, Goldsmith} that the components of a homotopically essential link in $S^3$ do not
bound disjoint maps of disks in $D^4$. The $2$-handle (a thickening of the disk) may be considered a trivial example of a Bing cell, and the content of theorem \ref{homotopically trivial} is the more general fact that the components of a homotopically essential link do not bound disjoint maps of Bing cells into the $4$-ball.

One may generalize further and consider surfaces (including those of higher genus) in $4$-manifolds with patches replaced by Bing cells. This leads to an interesting theory
that shares some of the features of both homology and homotopy. An important open question for applications to the A-B slice problem
(see section \ref{applications sec}) is to what extent this theory satisfies Alexander duality. This question will be pursued in a subsequent paper.

We will next summarize the ideas underlying the definition of the generalized Milnor group and the proof of theorem
\ref{homotopically trivial}. An outline of the argument here will be given to prove that the components $({\gamma}', {\gamma}'')$ of the Hopf link
do not bound disjoint maps into $D^4$ of two copies $C', C''$ of the Bing cell shown in figure \ref{cell figure}, where the maps are assumed to satisfy the disjointness requirements
discussed above. This example exhibits the main features of the general argument; a complete proof of theorem \ref{homotopically trivial} is given in section \ref{homotopy}.

Suppose to the contrary that there exist $C', C''$ whose attaching curves ${\gamma}', {\gamma}''$ form the Hopf link in $S^3=\partial D^4$.
Consider meridians $m'_i$ to
the $2$-handles $\{h'_i\}$ of $C'$ and meridians $m''_i$ to the $2$-handles $\{ h''_i\}$ of $C''$ in $D^4$.
By a meridian we mean a based loop in the complement representing the same homology class as a fiber
of the normal circle bundle over a $2$-handle. These meridians normally generate the fundamental group
${\pi}:={\pi}_1(D^4\smallsetminus (C'\cup C''))$. Consider the quotient $\overline M {\pi}$
of ${\pi}$ by the relations $[(m'_i)^x,( m'_j)^y]$ and $[(m''_i)^x, (m''_j)^y]$, where
$\{i,j\}\neq \{1,2\}, \{3,4\}$ and the conjugating elements $x,y$ range over all elements of ${\pi}$.
These relations arise from the double points between
the various $2$-handles in $D^4$ (they are an algebraic manifestation of the Clifford tori linking the double points, see section
\ref{surfaces in four ball} for further details).
The fact that the pairs $\{1,2\}, \{3,4\}$ are excluded reflects the fact that these
handles are required to be disjoint in the $4$-ball. The fact that the conjugates of
$m'_i, m''_j$ for any $i,j$ do not commute is due to the fact that $C', C''$ are disjoint.
This is a generalization of the Milnor group $M{\pi}$  introduced in \cite{Milnor1}, which in this context is defined
as the quotient of ${\pi}$ by the relations  $[(m'_i)^x,( m'_i)^y], [(m''_i)^x,( m''_i)^y], i=1,\ldots, 4$, corresponding to
{\em self}-intersections of the handles. These relations are included among the defining relations of $\overline M{\pi}$.

Let $m', m''$ denote the meridians to the components of the Hopf link in $\partial D^4$,
formed by the attaching curves ${\gamma}', {\gamma}''$. The relation $[m', m'']=1$ holds in ${\pi}_1(S^3\smallsetminus
({\gamma}'\cup {\gamma}''))$, therefore it also holds in ${\pi}_1(D^4\smallsetminus (C'\cup C''))$. On the other hand,
sliding the meridian along arcs in the pair of pants in figure \ref{cell figure} it is easy to see that the equalities
$m'=[m'_1,m'_2]=[m'_3,m'_4], m''=[m''_1,m''_2]=[m''_3,m''_4]$ hold in $\overline M {\pi}$. Carefully analyzing
the second homology of the complement $D^4\smallsetminus (C'\cup C'')$ one shows that there are no other
relations in $\overline M {\pi}$. Then the Magnus expansion can be used to prove that the
commutator $[m', m'']$ is in fact non-trivial in $\overline M{\pi}$: phrased differently, there are
``not enough'' relations in ${\pi}_1(D^4\smallsetminus (C'\cup C''))$ to imply the relation $[m', m'']=1$
which holds in ${\pi}_1(S^3\smallsetminus ({\gamma}'\cup {\gamma}''))$.
This contradiction concludes the outline of the proof that
the components of the Hopf link do not bound disjoint Bing cells of height $1$ in $D^4$. A complete proof of theorem \ref{homotopically trivial}
in section \ref{homotopy} relies on a detailed analysis of the relations in the generalized Milnor group of the complement
of Bing cells (of an arbitrary height) in $D^4$, which forms a central technical part of the paper.

\subsection{Applications: the A-B slice problem, topological arbiters} \label{applications sec}

The ideas introduced in this paper have been used to prove a number of results in
$4$-manifold topology. We will now summarize these applications.

The {\em A-B slice problem} is a reformulation of the $4$-dimensional topological surgery conjecture,
introduced by Freedman \cite{Freedman2} and further developed by Freedman-Lin \cite{FL}.
In this problem one considers decompositions of the $4$-ball, $D^4=A\cup B$, which extend the standard genus $1$
Heegaard decomposition of $S^3=\partial D^4$. The attaching curves ${\alpha}\subset \partial A, {\beta}\subset \partial B$
form the Hopf link in $\partial D^4$. The problem is then to find out whether there exist
decompositions $D^4=A_i\cup B_i$ and disjoint embeddings of the submanifolds $\{ A_i, B_i\}$ into $D^4$, so that
the attaching curves $\{ {\alpha}_i, {\beta}_i\}$ form a specified link (and its parallel copy) in $S^3$.
The central example of a link in question is the Borromean rings, see \cite{Freedman2, FL, K2} for more details.

In \cite{FL} the authors
introduced a family of {\em model decompositions} which appear to approximate,
in a certain algebraic sense, an arbitrary decomposition $D^4=A\cup B$.
In \cite{K2} the author showed how the idea of link groups of $4$-manifolds and theorem
\ref{homotopically trivial} may be used
to formulate an obstruction for the family of model decompositions:

\begin{theorem} \label{main AB theorem} \cite{K2} \sl
Let $L$ be the Borromean rings, or more generally any homotopically essential link
in $S^3$. Then $L$ is not $A$-$B$ slice where each decomposition $D^4=A_i\cup B_i$
is a model decomposition.
\end{theorem}

The proof is based on the observation that given a model decomposition, precisely one of the following two possibilities hold:
either ${\alpha}$ bounds a Bing cell in $A$ or ${\beta}$ bounds a Bing cell in $B$. Note that Bing cells of an arbitrary height,
not just height $1$, are needed to prove this result. Phrased in terms of link groups, either
${\alpha}=1\in {\lambda}(A)$ or ${\beta}=1\in{\lambda}(B)$. The proof then follows as a consequence of
theorem \ref{homotopically trivial}. This proof unified and generalized the previously known partial obstructions
in the $A$-$B$ slice program. The idea based on Bing cells is likely to be central in a solution to this problem.

The notion of a {\em robust} $4$-manifold is useful in putting theorems \ref{homotopically trivial}, \ref{main AB theorem}
in the context of link homotopy. Let $(M, {\gamma})$ be a pair
($4$-manifold, embedded curve in $\partial M$). The pair
$(M,{\gamma})$ is {\em robust} if whenever several copies $(M_i,
{\gamma}_i)$ are properly disjointly embedded in $(D^4, S^3)$, the
link formed by the curves $\{ {\gamma}_i\}$  in $S^3$ is
homotopically trivial. Therefore from the perspective of link homotopy theory, robust $4$-manifolds
act like disks (with self-intersections).
It follows from theorem \ref{homotopically trivial} that Bing cells are robust.
Theorem \ref{main AB theorem} may be rephrased as saying that
given a model decomposition $D^4=A\cup B$, precisely one of the two
parts $A$, $B$ is  robust.

We will now discuss an application of the results of this paper to the structure
of {\em topological arbiters}, established in \cite{FK}. Given an
$n$-dimensional manifold $W$, a topological arbiter associates a value $0$ or $1$ to
codimension zero submanifolds of $W$, subject to natural topological and duality axioms.
For example, there is a unique arbiter on ${\mathbb R}P^2$, which reports the location of the
essential $1$-cycle. The concept of a topolocial arbiter is rooted in
Poincar\'{e}-Lefschetz duality, indeed homology with field coefficients gives rise to arbiters
on projective spaces. A question addressed in \cite{FK} is the existence of arbiters not induced
by homology.

\begin{theorem} \label{uncountable} \cite{FK} \sl There exists an uncountable collection of local topological arbiters in dimension $4$.
\end{theorem}

Theorem \ref{homotopically trivial} is an important ingredient in the proof of this result.
A local arbiter is a version of a topological arbiter defined on the ball. It is defined
on codimension zero submanifolds of $D^4$ which meet $\partial D^4$ in a neighborhood of an unknotted
circle, and duality in this case is modeled on Alexander duality for homology.
Note that homology with various
field coefficients can be used to construct only a countable collection of arbiters (in any dimension).
Theorem \ref{uncountable} contrasts with
the situation in dimension $2$ where there is a unique local arbiter, and it is induced by homology.
A classification of topological arbiters (in dimensions other than $2$) remains an open problem; the tools
developed in this paper are an important ingredient in analyzing arbiters on $D^4$.

\subsection{Outline of the paper}
Sections \ref{nilpotent presentation}, \ref{facts} summarize the
relevant background material on presentations of nilpotent quotients,
Milnor's theory of link homotopy and related results on surfaces
in $4$-space. Section \ref{definitions section} introduces Bing
cells and link groups ${\lambda}(M)$, and gives examples of
$4$-manifolds with ${\pi}_1\neq{\lambda}$. Sections \ref{invariants}
-- \ref{homotopy} concern embeddings of Bing cells in $4$-space,
$(C,{\gamma})\hookrightarrow (D^4,S^3)$. More specifically, the
fundamental group of the complement, ${\pi}_1(D^4\smallsetminus C)$,
is analyzed in section \ref{invariants}.  Section \ref{generalized
Milnor} develops a generalization of the Milnor group in the context
of Bing cells. It is used, in particular, to define an algebraic
grading of Bing cells and link groups ${\lambda}_{i,j}$.
Sections \ref{representations}, \ref{invariant phi}, \ref{homotopy}
define an obstruction to embeddability of a collection of Bing
cells in $D^4$ with a prescribed boundary, given by a link in $S^3$.

\section{A presentation of nilpotent quotients}
\label{nilpotent presentation}

The purpose of this section is to describe a presentation of the
quotients ${\pi}/{\pi}^q$ of a group ${\pi}$ by the terms of its
lower central series, in terms of generators of the first and second
homology of ${\pi}$. This technique is well-known (see also
\cite{K}), and it will be used often throughout the paper. The lower central
series of a group $\pi$ is defined inductively by ${\pi}^1={\pi}$,
${\pi}^2=[{\pi},{\pi}], \ldots,{\pi}^q=[{\pi},{\pi}^{q-1}]$.

To state the lemma, fix a group $\pi$ and suppose that
$H_1(\pi;{\mathbb{Z}})$ is generated by $g_1,\ldots,g_n$,
$H_2(\pi;{\mathbb{Z}})$ is generated by $r_1,\ldots,r_m$, and let
$q\geq 2$ be an integer. Then the result of lemma~\ref{presentation}
is that, roughly speaking, $g_1,\ldots,g_n$ and $r_1,\ldots,r_m$ provide a
set of generators and relations respectively in a presentation of
$\pi/\pi^q$. To make this precise, consider the quotient homomorhism
$\alpha\co \pi/\pi^q\longrightarrow \pi/[\pi,\pi]$ and let $\hat
g_i\in \pi/\pi^q$ denote some preimage of $g_i$ under $\alpha$,
$i=1,\ldots,n$. It is a standard fact in nilpotent group theory
\cite{Warfield} that $\hat g_1,\ldots, \hat g_n$ generate
$\pi/\pi^q$.

Let $W\longrightarrow K(\pi,1)$ be a map from the wedge of $n$
circles $W$, inducing an epimorphism $\beta\co
\pi_1(W)\longrightarrow \pi/\pi^q$ and mapping the $i$-th free
generator of $\pi_1(W)$ to $\hat g_i$. Let $f_j\co
\Sigma_j\longrightarrow K(\pi,1)$ be a map of a surface $\Sigma_j$,
representing the generator $r_j$ of $H_2(K(\pi,1))\cong H_2(\pi)$,
$j=1,\ldots, m$. We assume here that each space has a fixed
basepoint, and all maps preserve them. The standard basis of
$H_1(\Sigma_j)$ pulls back via $\beta$ to some elements in
$\pi_1(W)$; let $\hat r_j\in \pi_1(W)$ be a lift via $\beta$ of the
attaching map of the $2$-cell of $\Sigma_j$. (In particular, if
$\Sigma_j$ is a $2$-sphere then the corresponding word $\hat r_j$ is
trivial.)

\smallskip

\begin{lemma} \label{presentation} \sl
Suppose $H_1(\pi; {\mathbb{Z}})$ is generated by $g_1,\ldots,g_n$,
and $H_2(\pi; {\mathbb{Z}})$ is generated by $r_1,\ldots,r_m$. Then
in the notations as above,
\[ \pi/\pi^q\cong \; <\!\hat g_1,\ldots,\hat g_n\: |\: \,
\hat r_1,\ldots,\hat r_m, (F_{\hat g_1,\ldots,\hat g_n})^q\!> \]

\noindent where $F_{\hat g_1,\ldots,\hat g_n}$ denotes the free
group on generators $\hat g_1,\ldots,\hat g_n$.
\end{lemma}

\smallskip

To prove the lemma we need a refinement of Stallings theorem
\cite{Stallings} due to Dwyer. Given a space $X$, the {\em Dwyer's
subspace} ${\phi}_k(X)\subset H_2(X;{\mathbb{Z}})$ is defined as the
kernel of the composition
$$ H_2(X)\longrightarrow H_2(K({\pi}_1 X,1))=H_2({\pi}_1 X)\longrightarrow
H_2({\pi}_1 (X)/{\pi}_1(X)^{k-1}). $$

\begin{theorem} \label{Dwyersthm}
{\rm \cite{Dwyer}} \sl Let $k$ be a positive integer and let $f\co
X\longrightarrow Y$ be a map inducing an isomorphism on $H_1(\, .
\,; {\mathbb{Z}})$ and mapping $H_2(X)/\phi_k(X)$ onto
$H_2(Y)/\phi_k(Y)$. Then $f$ induces an isomorphism
$\pi_1(X)/(\pi_1(X))^k$ $\cong\pi_1(Y)/(\pi_1(Y))^k$.
\end{theorem}

\smallskip

\noindent {\it Proof of lemma \ref{presentation}.} Let $X$ be the
$2$-complex obtained from $W$ by attaching $m$ two-cells along the
words $\hat r_1,\ldots \hat r_m$. The composition $W\longrightarrow
K(\pi,1)\longrightarrow K({\pi}/\pi^q,1)$ extends to $X$, inducing
an isomorphism $H_1(X)\cong H_1(\pi)\cong H_1(\pi/\pi^q)$ and an
epimorhism on $H_2/\phi_q$. Now an application of Dwyer's
theorem~\ref{Dwyersthm} concludes the proof of Lemma
\ref{presentation}. \qed

\section{The Milnor group: links in $S^3$ and surfaces in $D^4$} \label{facts}

In this section we recall the relevant results on Milnor groups and
$\bar\mu$-invariants \cite{Milnor1}, \cite{Milnor2}. This material
is used to set up the definition of Bing cells in section
\ref{definitions section}. Sections \ref{invariants} --
\ref{homotopy} develop a generalization of the Milnor group and of
other aspects of the theory in the context of Bing cells in
$D^4$.

\subsection{Links in $\mathbf{S^3}$} \label{links section}
Let $L=(l_1,\ldots,l_n)$ be an oriented link in $S^3$, and consider
meridians $m_1,\ldots,m_n$ to the components of $L$. By a meridian $m_i$ we mean
a path ${\gamma}_i$ in $S^3\smallsetminus L$ from a basepoint to the boundary of a regular
neighborhood of the component $l_i$, followed by a circle (a fiber of the circle normal bundle over $l_i$)
linking $l_i$ once and then followed by ${\gamma}_i^{-1}$ back to the basepoint.
Observe that
$H_1(S^3\smallsetminus L)$ is generated by $m_1,\ldots, m_n$, and a
set of generators for $H_2(S^3\smallsetminus L)$ is provided by
$n-1$ tori: the boundary of a regular neighborhood of $n-1$
components of $L$. By lemma \ref{presentation},
$\pi_1(S^3\smallsetminus L)/(\pi_1(S^3\smallsetminus L))^q$ has a
presentation
\begin{equation} \label{link presentation}
<m_1,\ldots,m_n|[m_1,w_1],\ldots,[m_{n-1},w_{n-1}],(F_{m_1,\ldots,m_n})^q>,
\end{equation}

where $F_{m_1,\ldots,m_n}$ denotes the free group generated by
$m_1,\ldots,m_n$, and $w_j$ is a word in $m_1,\ldots, m_n$ representing the untwisted $j$-th
longitude of the link.
The Magnus expansion homomorphism $M\co
F_{m_1,\ldots,m_n}\longrightarrow{\mathbb{Z}}\{x_1,\ldots,x_n\}$
into the ring of formal non-commutative power series in the
indeterminates $x_1,\ldots,x_n$ is defined by
$$M(m_i)=1+x_i, \, M(m_i^{-1})=1-x_i+x_i^2\pm\ldots$$

for $i=1,\ldots,n$. Let
\[ M(w_j)=1+\Sigma\mu_{L}(I,j)x_I \]

\noindent be the expansion of $w_j$, where the summation is taken over all
multi-indices $I=(i_1,\ldots,i_k)$ with entries between $1$ and $n$,
and $x_I=x_{i_1}\cdot\ldots\cdot x_{i_k}$, $k>0$. This expansion
defines for each such multi-index $I$ the integer $\mu_{L}(I,j)$. Let
$\Delta_{L}(i_1,\ldots,i_k)$ denote the greatest common divisor of
$\mu_{L}(j_1,\ldots,j_s)$ where $j_1,$ $\ldots,j_s$, $2\leq s\leq
k-1$ range over all sequences obtained by cancelling at least
one of the indices $i_1,\ldots,i_k$ and permuting the remaining
indices cyclically.

Let $\bar\mu_{L}(I)$ denote the residue class of $\mu_{L}(I)$ modulo
$\Delta_L(I)$. Analyzing the indeterminacy caused by the relations
in the presentation (\ref{link presentation}), one sees that for
each multi-index $I$ of length $|I|\leq q$ the residue class
$\bar\mu_{L}(I)$ is an {\it isotopy invariant} of the link $L$,
where $\bar\mu_{L}(I)$ is defined using the quotient
$\pi_1(S^3\smallsetminus L)/(\pi_1(S^3\smallsetminus L))^q$. In
particular, {\em the first non-vanishing coefficients ${\mu}_L(I)$
are well-defined.} (By first non--vanishing coefficients we mean
${\mu}_L(I)$ such that ${\mu}_L(J)=0$ for all proper subsets
$J\subset I$.)

\subsection{Link homotopy and Milnor groups.} \label{link-homotopy}

Two $n$-component links $L$ and $L'$ in $S^3$ are said to be {\em
link-homotopic} if they are connected by a 1-parameter family of
immersions such that different components stay disjoint at all
times. $L$ is said to be {\em homotopically trivial} if it is
link-homotopic to the unlink. $L$ is {\em almost homotopically
trivial} if each proper sublink of $L$ is homotopically trivial.

For a group $\pi$ normally generated by $g_1,\ldots,g_k$ its {\em
Milnor group} (with respect to $g_1,\ldots,g_k$) $M\pi$ is defined
to be the quotient of $\pi$ by the normal subgroup \begin{equation}
\label{Milnor relation} \ll [g_i,g_i^h]: 1\leq i\leq k, \; \;
h\in\pi\gg.\end{equation}

$M\pi$ is nilpotent of class $\leq k+1$, in particular it is a
quotient of $\pi/(\pi)^{k+1}$, and is generated by the quotient
images of $g_1,\ldots, g_k$. The Milnor group $M(L)$ of a link $L$
is defined to be $M\pi_1(S^3\smallsetminus L)$ with respect to its
meridians $m_i$.

Milnor showed in \cite{Milnor1} that the Magnus expansion induces a
well defined injective homomorphism $MM\co\!
M(F_{m_1,\ldots,m_k})\longrightarrow R(x_1,\ldots,x_k)$ into the
ring $R(x_1,\ldots,x_k)$ which is the quotient of
${\mathbb{Z}}\{x_1,\ldots,x_k\}$ by the ideal generated by monomials
$x_{i_1}\!\!\cdots\!x_{i_r}$ with some index occuring at least
twice. Indeed, every term in the Magnus expansion of each defining
Milnor relation (\ref{Milnor relation}) has repeating variables. Let
$\overline{w}_n\in MF_{m_1,\ldots,m_{n-1}}$ be a word representing
$l_n$ in $M\pi_1(S^3\smallsetminus (l_1\cup\ldots\cup l_{n-1}))$.
Then $\bar\mu$-invariants of $L$ with non-repeating coefficients may
also be defined by the equation
\[ MM(\overline{w}_n)=1+\Sigma\mu_{L}(I,n)x_{I} \]

\noindent where summation is over all multi-indices $I$ with
non-repeating entries between $1$ and $n\!-\!1$, and
$\bar\mu_{L}(I,n)$ is the residue class of $\mu_{L}(I,n)$ modulo the
indeterminacy $\Delta_{L}(I,n)$, defined above.

The Milnor group of $L$ is the largest common quotient of the
fundamental groups of all links link-homotopic to $L$, hence {\em if
$L$ and $L'$ are link homotopic then their Milnor groups are
isomorphic.} The next result gives an algebraic criterion for a link
to be null-homotopic.

\smallskip

\begin{lemma} \label{triviallink} {\rm \cite{Milnor1, Giffen, Goldsmith}} \sl
For an $n$-component link $L$, the following conditions are
equivalent:

    $(i)$ $L$ is homotopically trivial,\\
    $(ii)$ the components of $L$ bound disjoint immersed disks in $D^4$,\\
    $(iii)$ $M(L)\cong M(F_{m_1,\ldots,m_n})$ with the isomorphism carrying
    a meridian to $l_i$ to the generator $m_i$ of the free group,\\
    $(iv)$ all $\bar\mu$-invariants of $L$ with non-repeating coefficients vanish.
\end{lemma}

\smallskip

It follows from Lemma \ref{triviallink} that $L$ is almost
homotopically trivial if and only if all its $\bar\mu$-invariants
with non-repeating coefficients of length less than $n$ vanish. In
particular, if $L$ is almost homotopically trivial then its
$\bar\mu$-invariants with non-repeating coefficients of length $n$
are well-defined integers.

\subsection{The link composition lemma}

We will now recall the link composition lemma \cite{FL} (see also
\cite{KT}). The result on Bing cells proved in section
\ref{homotopy} contains this theorem as a special case. Given a link
$\widehat L =(l_1,\ldots, l_{k+1})$ in $S^3$ and a link
$Q=(q_1,\ldots, q_m)$ in the solid torus $S^1\times D^2$, their
``composition'' is obtained by replacing the last component of
$\widehat L$ with $Q$. More precisely, it is defined as
$C=(c_1,\ldots,c_{k+m}):=(l_1,\ldots,l_k,{\phi}(q_1),\ldots,{\phi}(q_m))$,
where ${\phi}\co S^1\times D^2\hookrightarrow S^3$ is a $0$-framed
embedding whose image is a tubular neighborhood of $l_{k+1}$. The
meridian $\{1\}\times\partial D^2$ of the solid torus will be
denoted by $\wedge$ and we set $\widehat Q:=Q\cup\wedge$.

\smallskip

\begin{theorem} \label{composition}  \sl If both $\widehat L$ and $\widehat Q$ are
homotopically essential in $S^3$ then so is their composition
$L\cup{\phi}(Q)$.
\end{theorem}

\subsection{Links in $\mathbf{S^1\times D^2}$} \label{links in torus}

Let $L$ be a link in $S^1\times D^2$. As above denote by $\wedge$ the
meridian $\{ p\}\times \partial D^2$ and set $\widehat L=L\cup
\wedge$. Consider $\widehat L$ as a link in $S^3$, using a standard
embedding $S^1\times D^2\subset S^3$. Links in the solid torus will
be used as attaching regions for $2$-handles in the definition of
Bing cells (see next section), and we need to specify the class
of links necessary for the definition. Let ${\wedge}'$ denote
another meridian $\{ q\}\times\partial D^2$, $p\neq q$.

\begin{definition} \label{admissible links} (Links used in the definition of Bing cells)
A link $L=(l_1\ldots,l_n)\subset S^1\times D^2$ is
{\em essential and (almost trivial)$^+$} if $\widehat L$ is homotopically
essential, and for each $i$, $(L\smallsetminus
l_i)\cup{\wedge}\cup{\wedge}'$ is homotopically trivial.
\end{definition}

An important example is given by $L=$Bing double of the core
$S^1\times\{ 0\}$ (see figure \ref{fig:surface Bing double}), and more generally by $L=$iterated Bing double of
the core. The fact that the (iterated) Bing doubles satisfy the conditions in definition \ref{admissible links}
follows from a computation of the $\bar\mu$-invariants of Bing doubles (cf. \cite{Milnor1, Cochran}), see the discussion below.
The definition also allows the trivial example: $L=$ core
of the solid torus.

The second condition is slightly stronger than just the requirement
that $\widehat L$ is almost homotopically trivial. We include it
since it is technically convenient for the proofs of the properties
of Bing cells. We need to reformulate the conditions on $L$ in
terms of $\bar\mu$ invariants. Consider the solid torus $S^1\times
D^2$ as the complement in $S^3$ of an unknotted circle and note that
$${\pi}_1((S^1\times D^2)\smallsetminus L)/ (( {\pi}_1(S^1\times D^2)\smallsetminus L))^q
\; \cong\; {\pi}_1(S^3\smallsetminus
(L\cup{\wedge}'))/{\pi}_1(S^3\smallsetminus (L\cup{\wedge}'))^q.$$

These groups are generated by the meridians $m_1,\ldots,m_n$ to the
components of $L$ and by the longitude $l=S^1\times \{ x\}$ of the
solid torus (respectively the meridian $\bar m$ to ${\wedge}'$ for the
second group.) Consider the free group $F_{m_1,\ldots,m_n,\bar m}$
mapping onto these groups, and the Magnus expansion
\begin{equation} \label{Magnus torus} M\co F_{m_1,\ldots,m_n,\bar m}\longrightarrow {\mathbb Z}
\{x_1,\ldots,x_n,y\}, \;\; M(m_i)=1+x_i,\, M(\bar m)=1+y.
\end{equation}

Let $W$ be a word representing $\wedge$ in the free group. Assuming
that $L$ satisfies the conditions in the definition above, observe
that all terms with {\em non-repeating variables} in the expansion
$M(W)$ are either of the form $x_{i_1}\cdots x_{i_n}$ or they
contain all variables $x_1,\ldots x_n$ and $y$. Since the link
$\widehat L$ is homotopically essential, renumbering the components
of $L$ if necessary, one can assume that the term ${\mu}\, x_1\cdots
x_n$ in the Magnus expansion $M(W)$ has the coefficient ${\mu}\neq
0$. It is important to note that there are no terms that contain $y$
and just a proper subset of the variables $x_1,\ldots,x_n$.

\subsection{Surfaces in $\mathbf{D^4}$.} \label{surfaces in four ball}

Let ${\Delta}=\cup {\Delta}_i$ be a collection of immersed disks in
$(D^4,\partial D^4)$. By Alexander duality, $H_1(D^4\smallsetminus
{\Delta})$ is generated by the meridians to the components of
$\Delta$, and $H_2(D^4\smallsetminus{\Delta})$ is generated by the
{\em Clifford tori} linking the double points of $\Delta$.

More precisely, a local model for the surfaces near a double point
is given by ${\mathbb R}^2\times\{ 0\}\cap \{ 0\}\times {\mathbb
R}^2\subset {\mathbb R}^2\times {\mathbb R}^2$. The Clifford torus
is the product of the unit circles $S^1\times S^1$. The linking
number of classes $a\in H_1({\Delta})$ and $b\in
H_2(D^4\smallsetminus{\Delta})$ may be computed as the intersection
number of ${\Sigma}\cdot b$ where $a=\partial {\Sigma}\subset D^4$.
$H_1({\Delta})$ is generated by the double point loops (loops in
$\Delta$ passing exactly once through a double point). It is clear
from the local model that the double point loops are paired up
${\delta}_{i,j}$ with the Clifford tori.

Suppose the disks ${\Delta}_i$ are disjoint, so all double points
are self-intersections. According to lemma \ref{presentation},
${\pi}_1(D^4\smallsetminus{\Delta})/({\pi}_1(D^4\smallsetminus{\Delta}))^q$
is generated by the meridians $m_1,\ldots, m_n$ to the components of
$\Delta$, and the relations (corresponding to the Clifford tori) are
all of the form $[(m_i)^f,(m_j)^g]=1$ for some $f,g$. In particular,
the Milnor group $M{\pi}_1(D^4\smallsetminus{\Delta})$ (with respect
to the meridian generators) is isomorphic to the free Milnor group
$MF_{m_1,\ldots,m_n}$.

This gives a useful perspective on the relation between $(i)$ and
$(ii)$ in lemma \ref{triviallink}: if a link $L$ is homotopically essential then
$M(L)$ is not isomorphic to the free Milnor group. This implies that
the components of $L$ do not bound disjoint maps ${\Delta}$ of disks
in $D^4$: otherwise the inclusion map $S^3\smallsetminus
L\longrightarrow D^4\smallsetminus {\Delta}$ would induce a
homomorphism $M(L)\longrightarrow MF_{m_1,\ldots,m_n}$, a
contradiction.

\section{Bing cells and link groups} \label{definitions section}

There are two kinds of Bing cells ($b$-cells) defined in this section. First we discuss abstract 
(``model'') Bing cells, then definition \ref{cell} introduces the notion of a Bing cell in a $4$-manifold $M$
which is a map of a model Bing cell into $M$ with only certain types of allowed singularities.

\begin{definition} \label{model cell}
A {\em model Bing cell ($b$-cell) of height $1$} is a smooth
$4$-manifold $C$ with boundary and with a specified attaching curve
${\gamma}\subset
\partial C$, defined as follows. Consider a planar surface $P$ with
$k+1$ boundary components ${\gamma}, {\alpha}_1,\ldots,{\alpha}_k$
($k\geq 0$), and set $\overline P=P\times D^2$. Let $L_1,\ldots,
L_k$ be a collection of links, $L_i\subset {\alpha}_i\times D^2$,
$i=1,\ldots,k$. We assume that for each $i$, $\widehat L_i$ is
essential and (almost trivial)$^+$, see definition \ref{admissible
links}. Then $C$ is obtained from $\overline P$ by attaching
zero-framed $2$-handles along the components of $L_1\cup\ldots\cup
L_k$.
\end{definition}

The surface $P$ (and its thickening $\overline P$) will be referred
to at the {\em body} of $C$, and the $2$-handles are the {\em
handles} of $C$.

A {\em model $b$-cell $C$ of height $h$}, $h>1$, is obtained from a $b$-cell of height $h-1$ by replacing its handles with $b$-cells of
height one. The {\em body} of $C$ consists of all (thickenings of)
its surface stages, except for the handles.

The most important example of links $L_i$ in the definition above is given by (iterated) Bing
doubles of the core ${\alpha}_i\times\{ 0\}$ of the solid torus ${\alpha}_i\times D^2$. These are the links
that appear in applications to the A-B slice problem \cite{K2} and to topological arbiters \cite{FK}.
Figure \ref{cell figure} in the introduction gives an example of a $b$-cell of height $1$: a schematic
picture and a precise description in terms of a Kirby diagram. Here
$P$ is a pair of pants, and each link $L_i$ is the Bing double of
the core of the solid torus ${\alpha}_i\times D^2$, $i=1,2$. The
reader is urged to keep this example in mind while reading the
paper: the theory already exhibits
many of its interesting features in this case.

{\bf Remarks.} 1. The standard $2$-handle $H=D^2\times D^2$ with
${\gamma}=\partial D^2\times\{ 0\}$ provides a trivial example of a $b$-cell (of any height) -- corresponding to the case $k=0$ in the
definition above. Alternatively, one gets the $2$-handle $H$ by
considering the links $L=$ cores of the corresponding solid tori.
Similarly, a $b$-cell of height $h$ also satisfies the definition
of a $b$-cell of any height $>h$.

2. One may assume that no body surface of $C$ above the first stage
is an annulus: suppose an annulus $A$ is present, $\partial
A={\gamma}_A\cup{\alpha}_A$.  Then $A$ may be used to deform the
attaching maps of handles or higher stages from ${\alpha}_A\times
D^2$ to ${\gamma}_A\times D^2$. This eliminates $A$ (and increases
the number of components of the link one stage below -- note that it
is still essential and (almost trivial)$^+$, see link composition
lemma \ref{composition} and also section \ref{homotopy}). So while
technically annuli are allowed by the definition, only planar
surfaces with $\geq 3$ boundary components above the first stage
contribute to the ``non-trivial'' increase of the height of $C$.

3. If the links $L$ defining $C$ have at least two components, then
$C$ is homotopy equivalent to the wedge of a collection of circles
and of a collection of $2$-spheres (one for each handle of $C$). Of
course if one considers $C$ up to homotopy then all relevant information (in definition \ref {admissible links})
about the attaching maps of the $2$-handles
is lost. This is the reason for the fact that link groups defined further below are a topological but not
in general a homotopy invariant of $4$-manifolds. Also note that in non-trivial examples of Bing cells $C$, $\gamma\neq 0\in H_1(C)$;
the link group ${\lambda}(M)$ and the first homology group $H_1(M)$ are not correlated.

4. In the definition above we used zero framed $2$-handles. In fact, in
light of definition \ref{cell} the framing is not going to be important for applications.

5. Recall the assumptions on each link $L$ in definition \ref{model
cell}: $(i)$ $\widehat L$ is homotopically essential, and $(ii)$
$L\cup {\wedge}\cup{\wedge}'$ is almost trivial. It is crucial for
the applications of $b$-cells that the link $L\cup{\wedge}$ is
essential -- this is made precise using the Magnus expansion
$M({\wedge})$, see section \ref{links in torus}. Therefore the basic
requirements on $L$ should be: $\widehat L$ is homotopically
essential and almost trivial. We made a slightly stronger
assumption: $L$ is (almost trivial)$^+$ since this makes the proofs
of the properties of $b$-cells technically easier. It is an
interesting question whether this extra condition may be removed in the proof of theorem
\ref{homotopically trivial}.

\subsection{The associated tree} \label{associated tree}
\begin{figure}[ht]
\vspace{.4cm} \centering
\includegraphics[height=2.5cm]{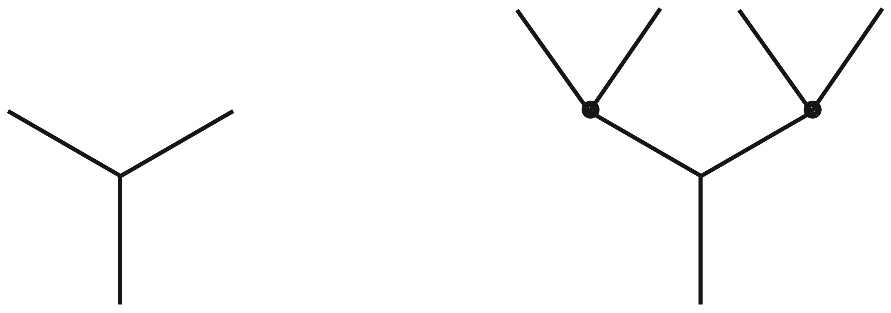}
    \put(-205,8){$T_P$}
    \put(-77,8){$T_C$}
{\Small
    \put(-97,65){$h_1$}
    \put(-69,65){$h_2$}
    \put(-30,65){$h_3$}
    \put(0,65){$h_4$}
    \put(-37,2){$\gamma$}}
    \vspace{.2cm} \caption{The trees $T_P, T_C$ corresponding to the Bing cell in figure \ref{cell figure}}
\label{tree example figure}
\end{figure}
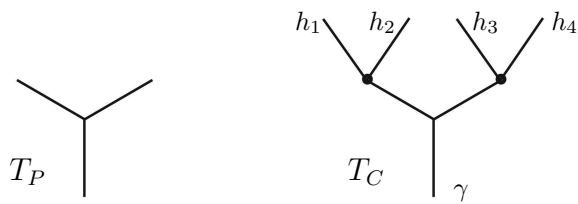

It is helpful to encode the
branching of a $b$-cell $C$ using an associated tree $T_C$ as follows. Define
$T_C$ inductively: suppose $C$ has height $1$. Then assign to the
body surface $P$ (say with $k+1$ boundary components) of $C$ the
cone $T_P$ on $k+1$ points. Consider the vertex corresponding to
the attaching circle ${\gamma}$ of $C$ as the root of $T_P$, and the
other $k$ vertices as the leaves of $T_P$. For each handle of $C$
attach an edge to the corresponding leaf of $T_P$, see figure \ref{tree example figure}.
The leaves of the resulting tree $T_C$ are in $1$-$1$ correspondence with the handles of
$C$.

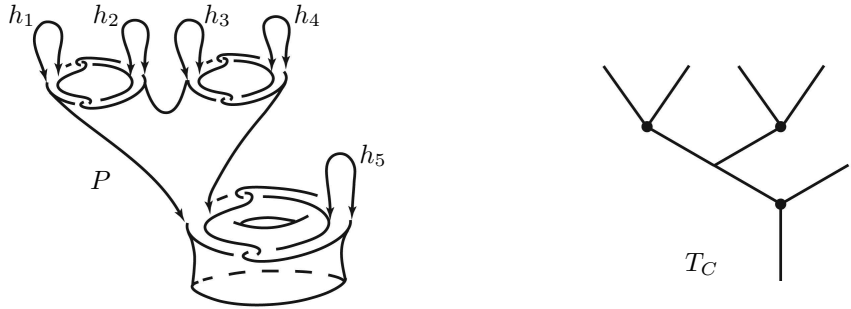
\begin{figure}[ht]
\centering
\includegraphics[height=4cm]{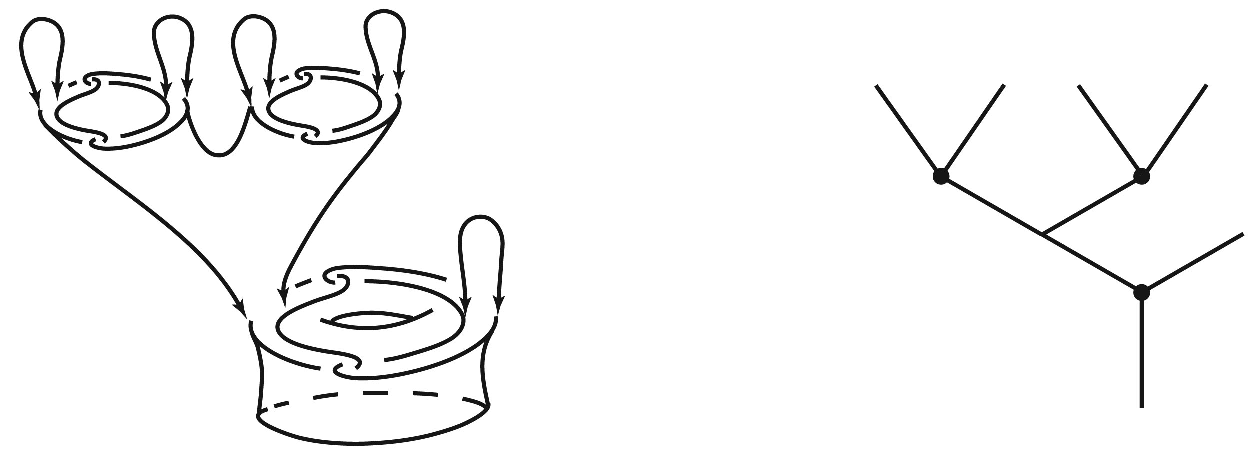}
{\Small
    \put(-65,15){$T_C$}
    \put(-321,109){$h_1$}
    \put(-289,109){$h_2$}
    \put(-247,109){$h_3$}
    \put(-213,109){$h_4$}
    \put(-188,56){$h_5$}
    \put(-290,46){$P$}}
    \vspace{.2cm} \caption{A schematic picture of a Bing cell $C$ of height $2$ and the
associated tree. Note that the bottom stage planar surface in $C$ is an annulus, giving rise to an
unmarked $2$-valent vertex which is not indicated in the bottom edge in $T_C$. }
\label{cell and tree}
\end{figure}

Suppose $C$ has height $h>1$, then it is obtained from a $b$-cell
$C'$ of height $h-1$ by replacing the handles of $C'$ with $b$-cells
$\{ C_i\}$ of height $1$. Assuming inductively that $T_{C'}$ is
defined, one gets $T_C$ by replacing the edges of $T_{C'}$
associated to the handles of $C'$ with the trees corresponding to
$\{ C_i\}$. Figure \ref{cell and tree} gives an example of a $b$-cell of height $2$ and its associated tree.

It is convenient to divide the vertices of $T_C$ into two types: the cone points
corresponding to body surfaces are {\em unmarked}; the rest
of the vertices are {\em marked} and are represented in figures by a wider dot.
Therefore the valence of an
unmarked vertex equals the number of boundary components of the
corresponding body surface. The marked vertices are in $1$-$1$
correspondence with the links $L$ defining $C$, and the valence of a
marked vertex is the number of components of $L$ plus $1$. It is
convenient to consider the $1$-valent vertices of $T_C$ (its root
and leaves, corresponding to the handles of $C$) as unmarked. This
terminology is useful in defining the maps of $b$-cells into $4$-manifolds below. The
height of a $b$-cell $C$ may be read off from $T_C$ as the maximal
number of marked vertices along a geodesic joining a leaf of $T_C$
to its root.

\subsection{Convention:} \label{preferred planar}
Recall from section \ref{links in torus} that for each link $L$ in the definition
of a Bing cell there is an ordering $l_1,\ldots, l_n$ of its
components so that the coefficient of the monomial $x_1\cdots x_n$ in the expansion
$M({\wedge})$ is non-trivial. Fix a specific planar embedding of $T_C$
reflecting this order, so that the clockwise ordering of the edges coincides with the ordering
$1,\ldots,n$ of the link components. This applies to marked vertices; there is a
flexibility in the planar embedding of the tree at its unmarked vertices.

\begin{definition} \label{cell}
A {\em Bing cell in a $4$-manifold $M$} is an embedding  $\overline C\subset M$, where
$C$ is a model Bing cell and $\overline C$ is the result of a finite number of
self-plumbings and plumbings among the handles and body surfaces of
$C$, subject to the following disjointness requirement:

$\bullet$ Let $A,B$ be either handles or body
stages of $C$, and let $a,b$ denote the corresponding vertices in $T_C$.
(For body surfaces this is the corresponding unmarked cone point,
for handles this is the associated leaf.) Consider the geodesic
joining $a,b$ in $T_C$, and look at its vertex $c$ closest to the
root of $T_C$. In other words, $c$ is the first common ancestor of
$a,b$. If $c$ is a marked vertex then no plumbings are allowed between $A$ and $B$.
\end{definition}

In particular, self-plumbings of any handle and body surface are
allowed. In the example shown in figure \ref{cell and tree}, the handle $h_1$ is
required to be disjoint from $h_2$, $h_3$ is disjoint from $h_4$;
$h_1$-$h_4$ and $P$ are disjoint from $h_5$. Abusing the notation, throughout the paper
we will denote Bing cells in $4$-manifolds by $C\subset M^4$, meaning the embedding of
the plumbed version $\overline C$ into $M$.

Note that a Bing cell $C$ is a thickening of a $2$-dimensional spine, and
in particular the solid tori ${\alpha}\times D^2$ which serve as the
attaching regions for higher stages are thickenings of circles.
From this perspective, given any map $f\co C\longrightarrow M^4$, it may be perturbed
so that all singularities are plumbings (thickenings of double points between handles and body surfaces),
and solid tori ${\alpha}\times D^2$ discussed above are embedded and disjoint from everything
else. The essential restriction in definition \ref{cell} is that handles and higher stages attached to different
components of each link $L_i$ defining a Bing cell are disjoint. It is straightforward to see that omitting this
restriction (i.e. allowing plumbing of arbitrary handles and body surfaces) would yield a trivial theory,
since any homotopically essential link may be unlinked by a suitable homotopy. On the other hand, there is no
disjointness requirement on handles/surfaces attached to different boundary components of a body surface.

\begin{definition} \label{link group}
Let $M$ be a $4$-manifold with a basepoint. Given $n\geq 1$, the $n$-th {\em link
group} ${\lambda}_n(M)$ is defined as $\{$based loops in
$M\}/{\sim}$. The equivalence relation $\sim$ on based loops in $M$ is defined as follows:
${\gamma}\sim {\gamma}'$ if there is
a based homotopy from ${\gamma}({\gamma}')^{-1}$ to a based loop
which bounds a Bing cell of height $n$ in $M$.
\end{definition}


\begin{proposition} \sl The relation ${\gamma}\sim{\gamma}'$ in
definition \ref{link group} is an equivalence relation, and moreover it is preserved by the product structure on loops.
\end{proposition}

{\em Proof.} Consider the first part of the statement, specifically the implication
$\gamma_1\sim\gamma_2, \gamma_2\sim\gamma_3\Rightarrow \gamma_1\sim\gamma_3$.
Assume $\gamma_1(\gamma_2)^{-1}$ is homotopic to a loop bounding a $b$-cell $C'$ and
$\gamma_2(\gamma_3)^{-1}$ is homotopic to a loop bounding $C''$, then
$\gamma_1(\gamma_3)^{-1}$ is homotopic to a loop bounding the wedge $(C_1,{\alpha}_1)\vee_p
(C_2,{\alpha}_2)$ of two $b$-cells of height $n$, where the
identification point $p$ is the base point. Using a boundary connected
sum of the bottom-stage surfaces, $C'\vee C''$ is converted into a $b$-cell of height $n$. Using isotopy, the attaching
regions of the form ${\alpha}\times D^2$ for higher-stage surfaces and handles of $C'$, $C''$ are made disjoint
from each other, since they are thickening of $1$-manifolds in $M^4$. The intersections
between arbitrary handles and body surfaces of $C'$ and those of $C''$ are allowed, since there is no
disjointness requirement on handles/surfaces attached to different boundary components of a body surface in definition \ref{cell}.

To prove the second part of the proposition, one needs to verify that if ${\gamma}_1\sim
{\gamma}'_1$ and ${\gamma}_2\sim {\gamma}_2'$ then
${\gamma}_1{\gamma}_2\sim{\gamma}_1'{\gamma}_2'$. This follows from
the equivalences ${\gamma}_1{\gamma}_2\sim{\gamma}_1'{\gamma}_2\sim{\gamma}_1'{\gamma}_2'$. \qed

\smallskip

{\bf Remarks.} 1. In light of remark 1 following definition \ref{model cell}, it follows that ${\pi}_1(M)$
surjects onto ${\lambda}_1(M)$.  Moreover, since a $b$-cell of
height $n$ satisfies the definition of a $b$-cell of height $n+1$, ${\lambda}_n(M)$
maps onto ${\lambda}_{n+1}(M)$. In section \ref{generalized Milnor}
we introduce an additional grading on $b$-cells, leading to a
two-parameter family of groups ${\lambda}_{i,j}(M)$.

2. Note that the definition of ${\lambda}_n(M)$ makes sense for a
manifold of any dimension (and in fact for any topological space),
but the theory is non-trivial only in dimension $4$. If $dim \, M\geq
5$ then the disjointness requirement is satisfied by general
position. If $dim\,  M<4$ then one doesn't expect it to hold due to the
dimension count.


It is easy to find examples of $4$-manifolds with ${\pi}_1\neq
{\lambda}_1$. Consider an example of $M^4$ with ${\pi}_1M\cong
{\mathbb Z}$ and ${\lambda}_1(M)=0$:

\begin{example} Consider $M=(S^1\times D^2\times I)\cup_L 2$-handles
where $L$ is the Bing double of the core of the solid torus $S^1\times
D^2\times\{1\}$, see figure \ref{fig:surface Bing double}. Clearly ${\pi}_1M\cong {\mathbb Z}$ and
${\lambda}_1(M)=0$. On the other hand, it is not difficult to see that $N=S^1\times D^3$ provides an
example where ${\lambda}_1(N)\cong{\pi}_1(N)\cong{\mathbb Z}$. See
lemma \ref{pi1l1} for a more detailed discussion.
\end{example}

\section{Bing cells in $4$-space} \label{invariants}

We begin the section by showing that any $b$-cell $(C,{\gamma})$ has
a realization in $(D^4,$ $\partial D^4)$. The main purpose of the
section is to analyze the fundamental group of the complement,
${\pi}_1(D^4\smallsetminus C)$. In particular, we will use the
technique presented in section \ref{nilpotent presentation} to find
a presentation of the nilpotent quotients ${\pi}_1(D^4\smallsetminus
C)/({\pi}_1(D^4\smallsetminus C))^q$. These results will be used in
sections \ref{generalized Milnor} -- \ref{invariant phi} to
formulate invariants which depend only on the underlying model
$b$-cell $C$ and not on its particular realization in the $4$-ball.
To fix the notation, recall that a model $b$-cell $(C,{\gamma})$ of
height $1$ is determined by the following data:

$\bullet$ the number of boundary components of the body surface $P$:
$\partial P={\gamma}\cup{\alpha}_1\cup\ldots\cup{\alpha}_k$,

$\bullet$ A collection of links $L_1\ldots,L_k$ where $L_i\subset
{\alpha}_i\times D^2$.



\begin{lemma} \label{embedding} \sl
Let $(C,{\gamma})$ be a model Bing cell. Then there is a realization
of $(C,{\gamma})\subset (D^4,\partial D^4)$.
\end{lemma}

Strictly speaking, the claim of the lemma is that there is an embedding
$(\overline C,{\gamma})\subset(D^4,\partial D^4)$ as in definition \ref{cell}. Abusing the notation we refer
to the plumbed version of the Bing cell as $C$ as well.

{\em Proof. of lemma \ref{embedding}.} The proof is inductive, starting with the base surface
of $C$ and moving up. Start with an unknotted circle ${\gamma}$ in
$S^3$ and let ${\gamma}\times D^2$ bound a $2$-handle $D^2\times
D^2$ in $D^4$. Puncture the core of the handle to get an embedding
of the first stage planar surface $P$. Note that for each $i$,
${\alpha}_i\times D^2$ bounds a (just removed) $2$-handle $H_i$ in
the complement of $P$ in $D^4$. The link $L_i\subset
{\alpha}_i\times D^2$ is homotopically trivial, so it bounds
disjoint immersed disks $\{ {\Delta}\}$ in $H_i$. The
self-intersections of handles and of body surfaces are allowed in
the definition of $b$-cells. If the height of $C$ is greater than
$1$, repeat the construction. (The disks $\Delta$ are converted into
second stage surfaces by puncturing them in the complement of double
points, etc.) \qed

\subsection{} \label{presentation section}
{\bf A presentation of $\mathbf{\mbox{\boldmath$\pi$}_1
(D^4\smallsetminus C)/(\mbox{\boldmath$\pi$}_1(D^4\smallsetminus
C))^q}$}.

For the remainder of this section fix $q\geq 2$. Given $(C,{\gamma})\subset (D^4,\partial D^4)$, let $m$ denote a
meridian to $\gamma$ in $S^3$. First assume $C$ has height one and
$P$ is a pair of pants, see figure \ref{cell figure} in the introduction.
Fix the notations: $$\partial
P={\gamma}\cup{\alpha}_1\cup{\alpha}_2,\,  \overline P=P\times D^2,\,
C=\overline P\cup_{L_1\cup L_2} 2{\rm -handles},$$ where $L_i\subset
{\alpha}_i\times D^2$ are links satisfying the conditions in
definition \ref{model cell}. Let $I_1$, $I_2$ be the index sets for
the components of $L_1,L_2$ respectively. ${\wedge}_1, {\wedge}_2$
will denote the meridional curves $\{p_i\}\times \partial D^2$ of
the solid tori ${\alpha}_i\times D^2$.

$H_1(D^4\smallsetminus C)$ is generated by $\mathcal{M}=\{ m_i\}$:
the meridians to the handles of $C$ (in the sense of section \ref{surfaces in four ball}). The index sets $I_1,I_2$ also
parametrize the handles of $C$, and to be specific, divide the set of
meridians $\mathcal{M}$ into two subsets: $\mathcal{M}_{I_1},
\mathcal{M}_{I_2}$. Denote by $F_{\mathcal{M}}=F_{M_{I_1},M_{I_2}}$
the free group generated by the elements of $\mathcal M$, and
consider the Magnus expansion $M$:

\begin{equation} \label{Magnusexpansion}
\pi_1(D^4\smallsetminus C)/(\pi_1(D^4\smallsetminus C))^q \;
\overset{p}{\longleftarrow}\; F_{\mathcal{M}}=F_{M_{I_1},M_{I_2}}\;
\overset{M}{\longrightarrow}\; {\mathbb Z}\{
X\}={\mathbb{Z}}\{X_{I_1},X_{I_2}\} \end{equation}

where as in section \ref{links section}, $M(m_i)=1+x_i$, $i\in I_1\cup I_2$. We need to fix a specific
word in $F_{\mathcal{M}}$ representing the meridian $m$.
Observe that ${\wedge}_1$ and $m$ cobound a cylinder in $D^4\smallsetminus C$:
the circle normal bundle of $P$ in $D^4$, restricted to a path in
$P$ joining two points in ${\alpha}_1$ and $\gamma$. Therefore
$m,{\wedge}_1$ are conjugate in ${\pi}_1(D^4\smallsetminus C)$.
Consider ${\wedge}_1$ in ${\pi}_1({\alpha}_1\times D^2\smallsetminus
L_1)$ and consider the commutative diagram of Magnus expansions, induced by the inclusion map
$i\co S^1\times D^2\smallsetminus L_1\subset D^4\smallsetminus C$:

\begin{equation} \label{comm diagram}
\xymatrix{
    {\pi}_1(S^1\times D^2\smallsetminus L_1)/({\pi}_1(S^1\times D^2\smallsetminus L_1))^q
    \ar[d]^{i_*} &
    F_{{\mathcal M}_1\cup \bar m}
    \ar[l] \ar[d]^{i_\sharp} \ar[r]^{M_1} &
    {\mathbb Z}\{ X_{I_1}, y\} \ar[d]
    \\
    {\pi}_1(D^4\smallsetminus C)/({\pi}_1(D^4\smallsetminus C))^q
    &
    \ar[l] F_{\mathcal M}
     \ar[r]^M &
    {\mathbb Z}\{ X\}
    }
\end{equation}

where $\bar m$ and $y$ are as in section \ref{links in torus} and specifically in (\ref{Magnus torus}).
The homomorphism $i_\sharp$ maps $m_j$ to $m_j$ for each $j\in {\mathcal M}_1$, and it maps
$\bar m$ to some fixed pullback of $i_*(\bar m)$ in $F_{\mathcal M}$.

Denote by $W_1$ some word representing ${\wedge}_1$ in the free
group $F_{\mathcal{M}_1\cup \bar m}$, then $W:= i_\sharp(W_1)$ represents $m$ in
$F_{\mathcal M}$. Recall (see section \ref{links in torus}) that
each term with non-repeating variable in the expansion $M_1(W_1)$ contains all of the variables
$x_1,\ldots,x_n$, where $I_1=\{ 1,\ldots, n\}$. According to the commutative diagram above,
this is also true for $M(W)$. It is important to remember (see last paragraph in section
\ref{links in torus}) that specifically $M(W)$ contains the non-trivial term ${\mu}\,  x_1\cdots x_n$.

Given an element $g\in \pi_1(D^4\smallsetminus
C)/(\pi_1(D^4\smallsetminus C))^q$, consider a word $w$ representing
it in $F_{\mathcal{M}}$. As in the classical case of Milnor's
invariants of links, discussed in section \ref{facts}, the
coefficients of the Magnus expansion $M(w)$ in general are not
well-defined invariants of $g$. This is due to the choice of the
meridians generating the group, and due to the fact that the kernel of the
surjection from $F_{\mathcal M}$ is non-trivial. In the present context, compared to the
classical situation, the kernel involves more relations in
$\pi_1(D^4\smallsetminus C)/(\pi_1(D^4\smallsetminus C))^q$
reflecting the topology of Bing cells.

According to lemma \ref{presentation}, to see the relations in
$\pi_1(D^4\smallsetminus C)/(\pi_1(D^4\smallsetminus C))^q$ we need
to analyze the generators of $H_2(D^4\smallsetminus C)$. By
Alexander duality, $H_2(D^4\smallsetminus C)\cong H_1(C,{\gamma})$.
Note that $H_1(C,{\gamma})$ is generated $(1)$ by double point loops
corresponding to the intersections among the handles and body
surfaces, subject to the disjointness requirement in definition
\ref{cell}, and $(2)$ by $H_1(P,{\gamma})$. (Here we assume the
non-trivial case: each link $L_i$ consists of at least two
components, so $C$ is homotopy equivalent to the wedge of two
circles with a collection of $2$-spheres, one $2$-sphere for each handle.) We
will divide the corresponding dual generators of
$H_2(D^4\smallsetminus C)$ into four types, $(R_1) - (R_4)$, and
analyze the resulting indeterminacy in the coefficients of the
Magnus expansion (\ref{Magnusexpansion}).

$\mathbf{(R_1)}$ \hspace{.3cm} Clifford tori for the
self-intersections of any handle $H_i$ of $C$, $i\in I_1\cup I_2$.

The corresponding relations are of the form $[(m_i)^f,(m_i)^g]=1$,
$i\in I_1\cup I_2$, $f,g\in{\pi}_1(D^4\smallsetminus C)$ (see
section \ref{surfaces in four ball}), and are familiar from the
study of link homotopy and the classical Milnor group (see
\ref{link-homotopy}). Pulling back the relations to $F_{\mathcal
M}$, consider the ideal $\mathcal{I}_1$ generated by their images in
${\mathbb Z}\{ X\}$. Observe that each term (besides $1$)
of any element in the ideal $\mathcal{I}_1$ has repeating variables.

More precisely, note that for any $a\in F_{\mathcal{M}}$ the Magnus
expansions $M(a^{-1}m_i a)$ and $M(a^{-1}(m_i)^{-1} a)$ are of the
form $1+$terms containing $x_i$ (where $M(m_i)=1+x_i$.) The
commutator $[(m_i)^f,(m_i)^g]$ is a product of $(m_i)^g$ conjugated
by $(m_i)^f$ and $(m_i^{-1})^g$, therefore
$M([(m_i)^f,(m_i)^g])=1+$terms containing at least two entries of
$x_i$. Hence the monomials with non-repeating variables are
invariant under multiplication by a conjugate of the relation
$(R_1)$.

According to definition \ref{cell}, any handle attached to $L_1$ can
intersect any handle attached to $L_2$. The corresponding generators
of $H_2$ are

$\mathbf{(R_2)}$ \hspace{.3cm} Clifford tori for the intersections
between the $2$-handles $H_{i_1}$ and $H_{i_2}$, \\
\hspace*{1.3cm} where $i_1\in I_1, i_2\in I_2$.

These tori give relations $[(m_{i_1})^f,(m_{i_2})^g]=1$. Each term
of any element in the ideal generated by the Magnus expansion of
these relations has both variables $x_{i_1}$ and $x_{i_2}$, where
$i_1\in I_1,i_2\in I_2$.

$\mathbf{(R_3)}$   \hspace{.3cm} Clifford tori for the intersections
of any handle $H_i$ with the body surface \\ \hspace*{1.3cm} $P$,
and Clifford tori for the self-intersections of $P$.

These generators of $H_2$ impose the relations of the form
$[m_i^f,m^g]$, and of the form $[m^f,m^g]$. Here $m_i$ is a meridian
to a handle $H_i$, $i\in I_1\cup I_2$ and $m$ is a meridian to $P$.
Recall from the discussion at the beginning of section
\ref{presentation section} that each term in the expansion $M(m)$
contains each of the variables $X_{I_1}$. If $i$ above is an element
of $I_1$ then all terms in the expansion of $[m_i^f,m^g]$ contain a
repeating variable (one of the $X_{I_1}$). If $i$ is an element of
$I_2$ then each term in the expansion of $[m_i^f,m^g]$ contains both
variables $x_{i_1}$ and $x_{i_2}$ for some $i_1\in I_1,i_2\in I_2$.
In either case, the indeterminacy has already appeared as a
result of relations $(R_1), (R_2)$.

There is another type $(R_4)$ of generators of $H_2(D^4\smallsetminus C)$,
Alexander dual to $H_1(P,\partial P\cap S^3)\cong {\mathbb Z}$.
Since we assumed each link $L_i$ has at least two components, the
meridian ${\wedge}_i=\{p_i\}\times \partial D^2$ of the solid torus
${\alpha}_i\times D^2$ bounds a surface $S_i$ in $({\alpha}_i\times
D^2)\smallsetminus L_i$. (Consider the disk $\{ p_i\}\times D^2$.
Since ${\wedge}_i$ has the trivial linking number with each
component of $L_i$, the disk may be converted into a surface
disjoint from the link.)

A geometric representative for this class of $H_2(D^4\smallsetminus
C)$ is given by the surface $S_1\, \cup\, $annulus$\,\cup\, S_2$.
Here the annulus is cobounded by ${\wedge}_1$ and ${\wedge}_2$, and
is the circle normal bundle of $P$ in $D^4$, restricted to a path in
$P$ joining two points in ${\alpha}_1, {\alpha}_2$. As above, denote
by $W_1, W_2$ some words in the free group representing
${\wedge}_1, {\wedge}_2$. Then the corresponding relation is

$\mathbf{(R_4)}$ \hspace{.3cm} $(W_1)^g \, (W_2)^{-1}=1$.

Now consider the {\em general height $=1$ case:} $\partial
P={\gamma}\cup{\alpha}_1\cup\ldots\cup{\alpha}_n$. The relations are
directly analogous to those described above; in particular there are
$n-1$ relations of type $(R_4)$: $(W_1)^{g_1}(W_2)^{-1}=1, \ldots,
(W_{n-1})^{g_{n-1}}(W_n)^{-1}=1$.

{\bf The general case} (height $\geq 1$). Denote by $\mathcal M$ the collection of meridians
$\{ m_i\}$ to the handles of $C$, and by $X$ a corresponding collection of
variables $\{ x_i\}$. The double points of $C$ occur as
intersections of handles and body surfaces, subject to the
disjointness assumption in definition \ref{cell}. More precisely,
the general relations of types $\mathbf{(R_1)-(R_3)}$ are
represented by the Clifford tori for self-intersections of each
handle and body surface of $C$, and for intersections of any two
handles and/or body surfaces, such that the first common ancestor of
the corresponding vertices in $T_C$ is an unmarked vertex. Recall that the
generators of $H_1(D^4\smallsetminus C)$, and also the variables $X$
are in $1$-$1$ correspondence with the handles of $C$ and also with
the leaves of $T_C$. The analysis directly analogous to the above implies
that {\em each term of any element in the ideal generated by the
Magnus expansions of the relations $(R_1)-(R_3)$ either contains
repeating variables, or it contains variables $x_{i}$ and $x_{j}$
whose first common ancestor in $T_C$ is unmarked.}

There is also a collection of relations $\mathbf{(R_4)}$ for the
body surfaces of $C$. Each generator of $H_1($body of
$C,{\gamma})$ contributes a relation of type $(W_1)^g(W_2)^{-1}$ as
above.

\section{The generalized Milnor group.} \label{generalized Milnor}

Starting with a Bing cell $(C,{\gamma})\subset(D^4,\partial D^4)$ we
will derive invariants of $(C,{\gamma})$ {\em independent of the
embedding} into $D^4$. This feature of the invariants is
particularly important for applications to the A-B slice problem
\cite{K2}. Recall from lemma \ref{triviallink}
that if a link $L\subset S^3$ is homotopically trivial, its components bound
disjoint immersed disks ${\Delta}$ in $D^4$, and the Milnor
group $M{\pi}_1(D^4\smallsetminus{\Delta})$ is isomorphic to the free Milnor group.
In particular (see section \ref{link-homotopy}) the
coefficients in the Magnus expansion
$M{\pi}_1(D^4\smallsetminus{\Delta})\longrightarrow R[X]$ are
well-defined. In our setting $M{\pi}_1(D^4\smallsetminus C)$ is {\em
not} the free Milnor group. The goal is to analyze the indeterminacy
and to extract useful invariants.

Recall the notation: we fix a collection $\mathcal{M}$ of meridians
$\{ m_i\}$ to the handles of $C$, one for each handle. Then the elements of
$\mathcal{M}$ generate any nilpotent quotient of
${\pi}_1(D^4\smallsetminus C)$.

\begin{definition}
The {\em generalized Milnor group} $GM(C)$ denotes
${\pi}_1(D^4\smallsetminus C)$ modulo the normal closure of all
elements of the form \begin{equation} \label{Milnor relations}
[m^f,m^g],\;\; {\rm and} \;\; [m_1^f, m_2^g],\;\; {\rm where}\;\;
f,g\in{\pi}_1(D^4\smallsetminus C), \;\, m,m_1,m_2\in\mathcal{M},\;
\; {\rm and}\end{equation}

the first common ancestor of $m_1,m_2$ is unmarked (see definition
\ref{cell}).
\end{definition}

In particular, $GM(C)$ is a quotient of the classical Milnor group
$M{\pi}_1(D^4\smallsetminus C)$ defined using the set $\mathcal{M}$
of normal generators. Consequently, $GM(C)$ is nilpotent, and so is
generated by the elements of $\mathcal{M}$.

For example, consider a realization in $D^4$ of the $b$-cell in
figure \ref{cell and tree} in section \ref{definitions section}. Denoting by $m_i$ a meridian to the handle $h_i$,
$i=1,\ldots,5$, the relations in the definition of $GM(C)$ are:
$$ [m_i^f,m_i^g]=1,\;\, i=1,\ldots, 5,\;\, [m_1^f,m_3^g]=
[m_1^f,m_4^g]= [m_2^f,m_3^g]= [m_2^f,m_4^g]=1,$$

where $f,g\in{\pi}_1(D^4\smallsetminus C)$.
The definition of $M(C)$ incorporates the relations $(R_1)-(R_3)$ in
${\pi}_1(D^4\smallsetminus C)$, discussed in the previous section.
In the classical Milnor's theory, the free Milnor group has a
well-defined representation into (the units of) the ring of
polynomials where the terms have non-repeating variables. In the
next section we describe the analogous representation for $GM(C)$.
In the present setup there is also an additional indeterminacy, due
to the relations $(R_4)$, and this is analyzed in section
\ref{invariant phi}. It is convenient to define, analogously to the
classical case, the free Milnor group:

\begin{definition} The {\em free} generalized Milnor group
$GM(F_{\mathcal{M}})$ is defined to be the free group
$F_{\mathcal{M}}$ modulo the relations of the form (\ref{Milnor relations}).
\end{definition}

It follows that $GM(C)$ is the quotient of $GM(F_{\mathcal{M}})$ by
the relations $(R_4)$.
Analogously to the classical case, $MC(G)$ has the
following property.

\begin{proposition} \sl
Given a model $b$-cell $C$, there exists a realization $\overline
C\subset D^4$ of $C$ such that ${\pi}_1(D^4\smallsetminus \overline
C)\cong GM(\overline C)$.
\end{proposition}

\smallskip

{\em Proof.} Consider any realization $C'\subset D^4$ of $C$.
$GM(C')$ is nilpotent and finitely generated, and is therefore
finitely presented. That is, $MC(G)$ is isomorphic to
${\pi}_1(D^4\smallsetminus C')$ modulo a finite number of relations
(\ref{Milnor relations}). It is a standard observation that these relations may be
introduced by finger moves yielding plumbings and self-plumbings of
$C'$ of the allowed type. This gives $\overline C$ satisfying the proposition. \qed

\subsection{Grading of Bing cells}
Given $(C,{\gamma})\subset (D^4,\partial D^4)$, let $m$ denote a
meridian to $\gamma$ in $S^3$. There is no relation, in general,
between the height of $C$ and how deep $m$ is in the lower central
series of ${\pi}_1(D^4\smallsetminus C)$, or of $GM(C)$. For
example, let $(C_1,{\gamma}_1)$ be a $b$-cell of height $k$ where
each link is the Bing double of the core of the corresponding solid
torus (and the body surfaces are arbitrary -- to be specific
consider pairs of pants.) Also consider $(C_2,{\gamma}_2)$ of height
$1$ where the body surface is an annulus and the link is the
$k$-iterated Bing double of the core. Then $C_1, C_2$ have different heights,
while for both $i=1,2$, $m_i$
is in the $2^k$-th term of the lower central series of
${\pi}_1(D^4\smallsetminus C_i)$.

Define the {\em nilpotency class} of $C$ to be the least $k$ such that the
$k$-th term of the lower central series $GM(C)^k$ is trivial.
Assuming that each link in the definition of $C$ has at least two
components, it is clear that the nilpotency class of a $b$-cell of
height $k$ is at least $k+1$. Refining definition \ref{link group},
consider ${\lambda}_{i,j}(M)=\{$based loops in $M\}$ modulo loops
bounding $b$-cells of height $i$ and having nilpotency class $j$.
There is a commutative diagram of surjections
\[
\xymatrix{
    \pi_1(M)
    \ar[r]
     &
    {\lambda}_{1,2}(M)
     \ar[r] &
    {\lambda}_{1,3}(M) \ar[d] \ar[r]
    &
    {\lambda}_{1,4}(M) \ar[d] \ar[r]
    &
    \ldots
    \\
    &
    &
    {\lambda}_{2,3}(M) \ar[r]
    &
    {\lambda}_{2,4}(M) \ar[r] \ar[d]
    &
    \ldots \\
    &
    & & \vdots
    }\]

\begin{lemma} \label{pi1l1} \sl
Let ${\pi},{\lambda}$ be finitely presented groups, where ${\pi}$ is
aspherical, and suppose ${\pi}$ maps onto $\lambda$. Then there are
$4$-manifolds $M$ with ${\pi}_1 (M)\cong {\pi}$ and
${\lambda}_{1,2}(M)\cong {\lambda}$.
\end{lemma}

\smallskip

{\em Proof.} Consider an aspherical $2$-complex $K$ with
${\pi}_1K\cong {\pi}$. Replacing the cells of $K$ by $0$-, $1$- and
$2$-handles, one gets a $4$-manifold $N$ with boundary. Observe that
${\pi}_1(N)\cong {\lambda}_{1,2}(N)$: suppose there is a loop
${\gamma}\subset N$ trivial in ${\lambda}_1(N)$ but not in
${\pi}_1(N)$. Then ${\gamma}$ is homotopic to a loop ${\gamma}'$
which bounds a $b$-cell $C$ of height $1$. Denote the body surface
of $C$ by $P$, $\partial
P={\gamma}'\cup{\alpha}_1\cup\ldots\cup{\alpha}_n$. It follows that
${\alpha}_i\neq 1\in {\pi}_1(N)$ for some $i$. The link $L_i\subset
{\alpha}_i\times D^2$ has two components. Consider the $2$-spheres
$S_1, S_2$ formed by the cores of the handles $H_1, H_2$ of $C$
attached to the components of $L_i$, capped off by the
null-homotopies of the components of $L$ in ${\alpha}_i\times D^2$.
Due to the assumptions on the link, and since the handles $H_1, H_2$
are disjoint, the intersection of $S_1$, $S_2$ is non-trivial in
${\mathbb Z}{\pi}_1(N)$. This is a contradiction with the
asphericity assumption.

Consider a collection of elements
${\alpha}=\{{\alpha}_1,\ldots,{\alpha}_k\}$ in ${\pi}_1K$ such that
the quotient of ${\pi}_1$ by the normal closure of ${\alpha}$ is
isomorphic to $\lambda$. Represent ${\alpha}$ by embedded
curves in $\partial N$, then the $4$-manifold $M=N\cup_{\alpha} 2$-handles,
where the handles are attached to Bing doubles of the cores of
${\alpha}_i\times D^2\subset \partial N$, satisfies the proposition. \qed

Examples of $4$-manifolds $M$ for which ${\lambda}_{i,j}(M)\neq
{\lambda}_{i,j+1}(M)$ are considered in \cite{K2}. It is an
interesting question whether there are manifolds for which vertical
maps are not isomorphisms either.

\section{Representations of $GM(C)$.} \label{representations}

The purpose of this section is to analyze the indeterminacy of the Magnus expansion
due to the relations in the generalized Milnor group $GM(C)$.
Consider a set $\mathcal{M}=\{ m\}$ of generators of
$H_1(D^4\smallsetminus C)$ provided by meridians to the handles
of $C$. The elements of $\mathcal{M}$ are in $1$-$1$ correspondence
with the leaves of the associated tree $T_C$, and are parametrized by multi-indices
$I=(i_1,j_1,\ldots,i_n,j_n)$ where $n$ is the height of $C$, the
indices $i_k$ correspond to the branching of the planar surface
stages, and the indices $j_l$ correspond to the components of the
attaching links $L$. Phrased in terms of the associated tree, the indices $i_k$ (respectively
$j_l$) correspond to branching of the tree $T_C$ at unmarked (respectively marked) vertices.

\begin{definition} \label{non-repeated}
Consider the set $X=\{ x\}$ whose elements are in $1$-$1$
correspondence with the elements of $\mathcal{M}$.  Let $R[C]$
denote the quotient of the free associative ring ${\mathbb Z}\{ X\}
$ generated by $X$ by the ideal generated by the monomials
$M=x_{I_1}\cdots x_{I_k}$  such that

$\bullet$ either $M$ contains repeating variables, or

$\bullet$ $M$ contains variables $x_I$, $x_{I'}$ whose first common
ancestor in $T_C$ is unmarked (compare with definition \ref{cell}).
\end{definition}

The second condition may be rephrased as follows: let $I=(i_1,\ldots, j_n)$, $I'=(i'_1,\ldots,j'_n)$
be two multi-indices as above.
Consider the first index where these sequences differ: if it is one
of the $j$'s then any monomial containing $x_I$, $x_{I'}$ is in the
ideal.

\begin{proposition} \label{proposition diagram} \sl
The Magnus expansion $F_{\mathcal{M}}{\longrightarrow} {\mathbb Z}\{X\}$ induces a well-defined homomorphism
$GM(F_{\mathcal{M}})\longrightarrow R[C]$, which
abusing the notation  is also denoted $M$:
\[
\xymatrix{
    F_{\mathcal{M}}
    \ar[r] \ar[d]
     &
    GM(F_{\mathcal{M}})
     \ar[r] \ar[d]^{M} &
    GM(C)
    \\
    {\mathbb Z}\{X\}
    \ar[r] &
    R[C]&
    }\]
\end{proposition}

{\em Proof.} The kernel of $F_{\mathcal{M}}\longrightarrow
GM(F_{\mathcal{M}})$ is normally generated by the relations
(\ref{Milnor relations}). Note that every every term (besides $1$)
in the expansion $M([m^f,m^g])$ has a repeating variable, $x$, corresponding to $m$.
Similarly, every term in the expansion $M([m_1^f, m_2^g])$ contains
both variables $x_1, x_2$. Therefore the expansion of each relation
is in the ideal defining $R[C]$. \qed

\begin{definition} \label{rep R}
Let $v$ be a vertex of $T_C$. Assign to it
an additive subgroup $\widetilde R_v\subset R[C]$ as follows.
Denote by $T_v$ the subtree of $T_C$ rooted at $v$, and let $X_v$
denote all variables corresponding to the leaves of $T_v$.

At each unmarked vertex of $T_v$ keep exactly one branch and erase
the rest. Denote the resulting subtrees rooted at $v$ by $\{
T^{\alpha}_v\}$. Then $\widetilde R_v$ is defined to be the span of
the monomials read off, clockwise starting at the left-most leaf, from all possible planar
embeddings of $T^{\alpha}_v$, for all $\alpha$.

Suppose $C$ has height $n$ and without loss of generality assume all branches have uniform
length (insert extra stages = annuli if necessary). Set $\widetilde
R_k(C)=\oplus _v \widetilde R_v\subset R[C]$, where the direct sum is
taken over all vertices $v$ at height $k$. Denote $$\widetilde
R(C): =\widetilde R_0(C)=\widetilde R_{{\rm root\; \, of\; \,}T_C}.$$
\end{definition}

\begin{figure}[ht]
\centering
\includegraphics[width=12cm]{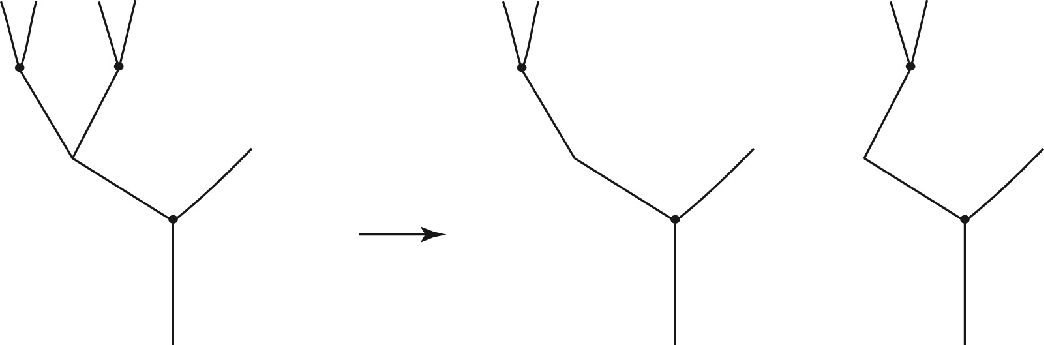} \hspace{1.2cm}
{\Small
    \put(-350,117){$x_1$}
    \put(-333,117){$x_2$}
    \put(-317,117){$x_3$}
    \put(-300,117){$x_4$}
    \put(-261,68){$x_5$}
    \put(-184,117){$x_1$}
    \put(-167,117){$x_2$}
    \put(-97,68){$x_5$}
    \put(-58,117){$x_3$}
    \put(-41,117){$x_4$}
    \put(-3,68){$x_5$}}
\caption{}
\label{cell figure trees}
\end{figure}

For example, consider the Bing cell in figure \ref{cell and tree}. Then there are two
subtrees entering the definition of $\widetilde R(C)$, shown in
figure \ref{cell figure trees}. There are a total of $8$ planar embeddings of these
subtrees, giving the monomials $\{ x_1x_2x_5,\, x_2x_1x_5,\,
x_5x_1x_2,\, x_5x_2x_1$, $x_3x_4x_5, \, x_4x_3x_5,\,
x_5x_3x_4,\, x_5x_4x_3\}.$ Some of the terms,
for example $x_1x_5x_2$, do not appear since they do not arise
from a subtree.

We will be interested in the subring $1+\widetilde R(C)$ of $R[C]$.
It follows from definition of $R[C]$ and from the assumptions in definition \ref{admissible links} on the links forming the Bing cell $C$
that all monomials in $\widetilde R(C)$ have ``maximal
length''. That is, if $X_I$ is a monomial in $\widetilde R(C)$ then for any
variable $x\in X$, inserting $x$ anywhere in $X_I$ gives a trivial
element of $R[C]$. Observe that the product in $1+\widetilde R(C)$ is given by
$$(1+\sum_I {\alpha}_I X_I)(1+\sum_I{\beta}_I X_I)=1+\sum_I
({\alpha}_I+{\beta}_I) X_I.$$

\begin{definition} \label{rep S}
For each vertex $v$ of $T_C$ consider the subring
$$S_v=1+\widetilde R_v+{\rm higher}\;\;{\rm order}\;\; {\rm
terms}$$

of $R[C]$. By higher order terms we mean all terms of the form
$$T=f_1\, x_1\, f_2\, x_2\, \ldots\, f_m\, x_m\, f_{m+1}$$
where the monomial $x_1\ldots
x_m$ (obtained from $T$ by deleting the $f$'s) is in $\widetilde
R_v$, and at least one of the monomials $f_1,\ldots,f_{m+1}\in
{\mathbb Z}\{ X\}$ is not equal to $1$. Similarly, set
$S_k=1+\widetilde R_k+$higher order terms. Observe that
$S_0(C)=1+\widetilde R(C)$: the monomials in $\widetilde R(C)$
already have maximal length, so there are no higher order terms.
\end{definition}

\begin{figure}[ht]
\vspace{.35cm}
\centering
\includegraphics[width=10cm]{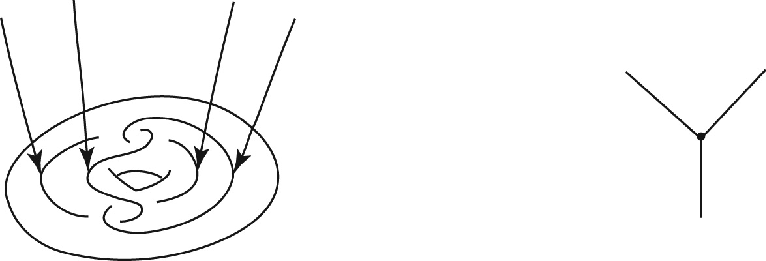} \hspace{1.2cm}
{\small    \put(-66,83){$T_{C_1}$}
    \put(-191,95){$C_2$}
    \put(-307,24){$C$}
    \put(-279,95){$C_1$}
    \put(-6,83){$T_{C_2}$}
    \put(-52,20){$T_C$}}
    \caption{Raising the height: step $1$.}
\label{step 1}
\end{figure}
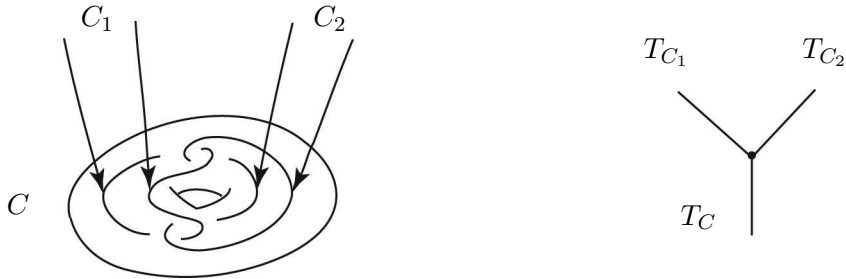

\subsection{} \label{inductive}
\begin{figure}[ht]
\vspace{.35cm}
\includegraphics[width=10cm]{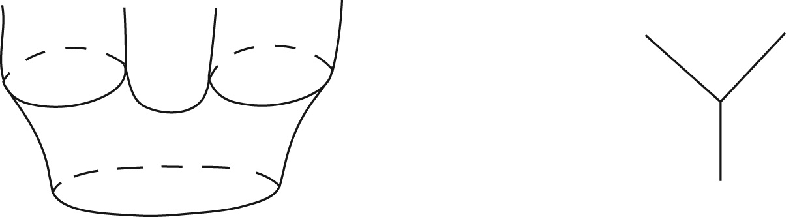} \hspace{1.2cm}
{    \put(-72,80){$T_{C_1}$}
    \put(-192,82){$C_2$}
    \put(-293,16){$C$}
    \put(-270,82){$C_1$}
    \put(-3,80){$T_{C_2}$}
    \put(-52,17){$T_C$}}
    \vspace{.2cm} \caption{Raising the height: step $2$.}
\label{step 2}
\end{figure}
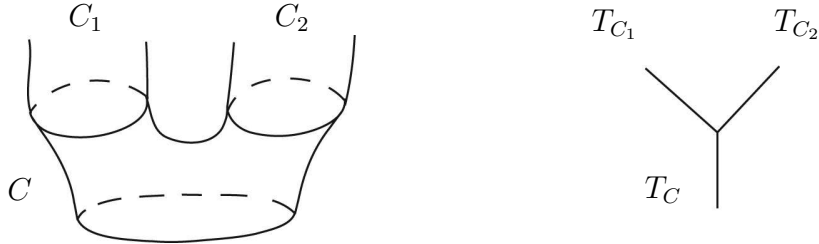

It is useful to note an inductive construction of the
representation $S_0(C)=1+\widetilde R(C)$. A Bing cell of height $k$
is assembled from a bottom stage planar surface $P$, $\partial
P={\gamma}\cup{\alpha}_1\cup\ldots\cup{\alpha}_n$, and Bing cells of
height $k-1$ attached to the components of the links $L_i$,
$L_i\subset {\alpha}_i\times D^2$. This assembly may be decomposed
into two steps. Step one (figure \ref{step 1}) corresponds to attaching
Bing cells of height $k-1$ to a single link $L$. Supposing for simplicity of notation that
$L$ consists of just two components, it follows from definition \ref{rep R}
that in this case $$\widetilde R(C)\cong (\widetilde
R(C_1)\otimes \widetilde R(C_2)) \oplus (\widetilde R(C_2)\otimes
\widetilde R(C_1)),$$ with the obvious generalization for links $L$ with more than two components.
Here the map $\widetilde R(C_i)\otimes
\widetilde R(C_j)\longrightarrow \widetilde R(C)$ is defined on
generators by $X_i\otimes X_j \longmapsto X_i\cdot X_j$, the product
of monomials. Step two (figure \ref{step 2}) combines the results of step
one which are attached to an arbitrary planar surface. In this case
$\widetilde R(C)\cong \widetilde R(C_1)\oplus \widetilde R(C_2)$.

\begin{lemma} \label{rep properties} \nl \sl
1. Let $m$ be a meridian to a body surface of $C$, and let $v$ be
the corresponding vertex in $T_C$. Then there exists a word $w\in
F_{\mathcal M}$ representing it so that $M(w)\in S_v$.

2. In particular, let $m_0$ denote a meridian to the bottom stage of
$C$ in $D^4$ (for example, a meridian to $\gamma$ in $S^3$.) Then
there exists a word $w_0$ representing it in $F_{\mathcal M}$ such
that $M(w_0)\in S_0(C)=1+\widetilde R(C)$.
\end{lemma}

{\em Proof.} The proof is inductive, moving from the handles down.
If $m$ is a meridian to a handle of $C$ then $M(v)=1+x$ and the
statement is obviously true. Suppose the statement holds for the
meridians to all body surfaces at height $k+1$, and let $m$ be a
meridian to a surface $P$ at height $k$. Note that the statement is
independent of a choice of the meridians: if one of the meridians is
replaced by a conjugate, the Magnus expansion would also satisfy the
condition. Denote, as usual, $\partial
P={\gamma}\cup{\alpha}_1\cup\ldots\cup{\alpha}_n$; the surfaces at
height $k+1$ are attached to $P\times D^2$ along the links $L_i$,
$L_i\subset {\alpha}_i\times D^2$. For each $i$, the meridian $m$ is
conjugate to the curve ${\wedge}_i$ (connected to the basepoint).
Therefore for the inductive step it suffices to consider only step
one of the height raising discussed above. In other words, one can
assume that $P$ is an annulus, and there is only one link $L\subset
{\alpha}\times D^2$.

Consider the map ${\pi}_1({\alpha}\times D^2\smallsetminus
L)\longrightarrow {\pi}_1(D^4\smallsetminus C)$. The map is obtained
by pushing ${\alpha}\times D^2\smallsetminus($a thickening of $L)$
slightly into the complement of $C$ in $D^4$. Let $L=(l_1,
\ldots,l_n)$; denote the corresponding Bing cells 
attached to them by $C_1,\ldots,C_n$, as in figure \ref{step 1}. To
distinguish them from the meridians to the handles of $C$, denote
the meridians to the components of $L$ in the solid torus by
$m'_1,\ldots,m'_n$, and let $z_1,\ldots,z_n$ be the corresponding
variables for the Magnus expansion. Denote the longitude of the
torus, $\{p\}\times \partial D^2$, by $l$, and the corresponding
variable by $y$.

The meridians $m'_j$ to the components of $L$ may be viewed as
meridians to the bottom surface stages of $C_j$. By the inductive
assumption, there are preimages $w_j$ of $i_*(m'_j)$ in $F_{\mathcal
M}$ such that the Magnus expansion $M(w_j)$, composed with the
projection to $R[C]$, is in $S_{v_j}=1+\widetilde R_{v_j}+$higher
order terms. In the following diagram, the map $\phi$ between the
free groups is defined on generators by taking the preimage
$w_j$ of $i_*(m'_j)$ in $F_{\mathcal M}$. Similarly,
${\phi}(l)$ is defined as a pullback of $i_*(l)$ in $F_{\mathcal M}$. Then ${\psi}(z_j)$ is
defined as $M({\phi}(m'_j))-1=M(w_j)-1$.

\[
\xymatrix{
    {\pi}_1({\alpha}\times D^2\smallsetminus L)/({\pi}_1({\alpha}
    \times D^2\smallsetminus L))^q
    \ar[d]^{i_*} &
    F_{m'_1,\ldots,m'_{n},l}
    \ar[l] \ar[d]^{\phi} \ar[r]^{M'} &
    {\mathbb Z}\{z_1,\ldots,z_n,y\} \ar[d]^{\psi}
    \\
    {\pi}_1(D^4\smallsetminus C)/({\pi}_1(D^4\smallsetminus C))^q
    &
    F_{\mathcal M}
    \ar[l] \ar[r]^M &
    {\mathbb Z}\{ X\} \ar[r]^{\pi} & R[C]
    }\]

Recall from the discussion preceding this lemma that $$ \widetilde
R_v\cong \oplus [\widetilde R_{v_{i_1}}\otimes\ldots\otimes
\widetilde R_{v_{i_n}}]$$

where the direct sum is taken over all permutations of $\{
1,\ldots,n\}$, and the inclusion $\widetilde
R_{v_{i_1}}\otimes\ldots\otimes \widetilde
R_{C_{i_n}}\longrightarrow \widetilde R_v$ is defined on the
additive generators by multiplication of the monomials.
Let $w$ be a word representing $\wedge$ (the meridian of the torus ${\alpha}\times D^2$)
in the free group $F_{m'_1,\ldots,m'_{n},l}$. We will
use the assumptions \ref{admissible links} on the links $L$ in the
definition of Bing cells \ref{model cell}. In particular, every
term with non-repeating variables in the expansion $M'(w)$ contains
each of the variables $z_1,\ldots,z_n$ (and in addition it may also
contain $y$.) The expansion $M({\phi}(w))$ is obtained from $M'(w)$
by replacing each $z_i$ and $y$ with ${\psi}(z_i)$, ${\psi}(y)$. The
proof is completed by the observation that
$${\pi}(\prod_{j=1}^n{\psi}(z_{i_j}))\; \; {\rm and} \; \; {\pi}[{\psi}(z_{i_1})\; \ldots\; {\psi}(z_{i_k})\;{\psi}(y)\;
{\psi}(z_{i_{k+1}})\;\ldots\; {\psi}(z_{i_n})]$$ are elements of
$S_v$, provided that for each $j$, ${\pi}({\psi}(z_{i_j}))\in S_{v_{i_j}}$.
The expansion $M'(w)\in
{\mathbb Z}\{z_1,\ldots,z_n,y\}$ may contain a {\em proper} subset
of the variables $\{ z_1,\ldots, z_n\}$, provided that at least one
of them, say $z_i$, is {\em repeated}. However by assumption
${\psi}(z_i)\in S_{v_i}$, so according to definition \ref{rep R}
every term of ${\psi}(z_i)$ contains all of the variables associated
to a subtree $T^{\alpha}_{v_i}$. Then to analyze ${\psi}(z_i)\cdots
{\psi}(z_i)$ consider the product of any two such terms. Either they
correspond to the same tree $T^{\alpha}$ and then the product
contains repeated variables and so is trivial in $R[C]$, or they
correspond to different subtrees $T^{\alpha}$, $T^{\beta}$, and then
the product is again trivial in $R[C]$, by the second condition in
definition \ref{non-repeated}. \qed


To define invariants of Bing cells in the next section, we need to
fix a more specific subspace of $\widetilde R_v$, for each $v$,
containing precisely the monomials with non-trivial
$\bar\mu$-invariants of the links $\widehat L$ in the definition of
$b$-cells (see definition \ref{model cell} and the discussion at the
end of section \ref{links in torus}.) The definition is similar to
that of $\widetilde R_v$ but it involves only a specific order of
the variables $X$.

\begin{definition} \label{rep Q}
Let $v$ be a vertex of $T_C$. Consider the subtrees $T^{\alpha}_v$
of $T_C$ whose root is $v$, as in definition \ref{rep R}. Then $Q_v$
is the additive subgroup of $R[C]$ spanned by the monomials read
off, clockwise, from the {\em fixed} planar embedding, defined in
\ref{preferred planar}, of $T^{\alpha}_v$, for all $\alpha$.
(Therefore $Q_v\subset \widetilde R_v$.) Set $Q_k(C)=\oplus _v
Q_v\subset R[C]$, where the summation is taken over all vertices $v$
at height $k$. Also denote $Q(C)=Q_0(C)=Q_r$ where $r$ is the root
of $T_C$.
\end{definition}

In the example in figures \ref{cell and tree}, \ref{cell figure trees}, $Q(C)$ is spanned by the monomials
$x_1x_2x_5,\; x_3x_4x_5$. (Compare with the computation of
$\widetilde R(C)$ in this example, following definition \ref{rep
R}.)

We will also use an alternative, inductive, description of $Q(C)$,
analogous to that of $\widetilde R(C)$ (see \ref{inductive}). For
each leaf $l$ of $T_C$, the corresponding $Q_l$ is the subgroup
($\cong{\mathbb Z}$) of $R[C]$ spanned by $x_l$. Suppose $Q_v$ is
defined for vertices of $T_C$ at height $>k$, and let $v$ be an
(unmarked) vertex at height $k$. Moving down the Bing cell from
height $k+1$ to height $k$ may be decomposed into steps, illustrated in figures \ref{step 1},
\ref{step 2}. The first step (corresponding to $P=$annulus) gives $Q\cong
Q_1\otimes Q_2$. The second step (figure 5) gives $Q\cong
Q_1\oplus Q_2$. To combine these two steps, denote $\partial
P={\gamma}\cup {\alpha}_1\cup\ldots\cup {\alpha}_n$; surfaces at
height $k+1$ are attached along the links $L_i\subset
{\alpha}_i\times D^2$. Let $I_i$ be the (ordered) index set for the
components of $L_i$. Then
\begin{equation} \label{tensor}
Q=\bigoplus_i\;  \bigotimes_{j\in I_i} \; \, Q_j
\end{equation}
{\em Remark.} The structure of $Q(C)$ may be read off from the tree
$T_C$ associated to $C$: the ``generators'' correspond to the leaves
of $T_C$; then form a tensor product for each marked vertex of the
tree and a direct sum for each unmarked vertex.

\subsection{The ring structure} \label{ring structure}

For each $v$, $S_v$ is a subring of $R[C]$. Consider $1+\widetilde
R_v$ as the quotient of $S_v$ by the ideal generated by the higher
order terms (see definition \ref{rep S}), and let ${p}_1\co
S_v\twoheadrightarrow 1+\widetilde R_v$ denote the projection.
Similarly, $1+Q_v$ is the quotient of $1+\widetilde R_v$ by the
ideal generated by all monomials which do not respect the fixed
order of the variables, ${p}_2\co 1+\widetilde R_v\twoheadrightarrow
1+Q_v$. The product in $1+\widetilde R_v$, $1+Q_v$ is given by
$$(1+\sum_I {\alpha}_I X_I)(1+\sum_I{\beta}_I X_I)=1+\sum_I
({\alpha}_I+{\beta}_I) X_I.$$ Let $m$ be a meridian to the bottom
stage of $C$, then by lemma \ref{rep properties} there exists
a word $w$ representing it in the free group whose Magnus expansion
$M(w)$ is an element of $S(C)$. Consider its image in $1+Q(C)$:
\begin{equation} \label{expansion}
{p}_2({p}_1(M(w)))=1+\sum_I {\alpha}_i X_I,
\end{equation}
where the summation is over all subtrees with a prescribed planar
embedding, as discussed above. The coefficients ${\alpha}_I$ are
well-defined with respect to the relations $(R_1)-(R_3)$ of section \ref{invariants}. (That is,
with respect to multiplying $w$ by a conjugate of one of the
relations $(R_1)-(R_3)$.) The next section introduces an invariant
well-defined with respect to $(R_4)$ as well.

\section{An invariant $\Phi$ of Bing cells.} \label{invariant phi}

The purpose of this section is to prove the following
statement. The main content is in the proof, which will be generalized from
knots to the setting of links in section \ref{homotopy}.

\begin{lemma} \label{invariant} \sl
Let $K$ be a knot in $S^3$, suppose $C$ is a Bing cell in $D^4$ bounded by $K$, and fix $q\geq 2$.
If $g$ is an element of ${\pi}_1(S^3\smallsetminus
K)$ whose image is non-trivial in $H_1(S^3\smallsetminus K)$, then  $i_*(g)\neq 1\in {\pi}_1(D^4\smallsetminus
C)/({\pi}_1(D^4\smallsetminus C))^q$. Here $i_*$ is the map induced by the inclusion
$i\co S^3\smallsetminus K\subset D^4\smallsetminus C$.
\end{lemma}

\subsection{} \label{special word}
{\em Notation.}
Given $g\in {\pi}_1(S^3\smallsetminus
{\gamma})/{\pi}_1(S^3\smallsetminus {\gamma})^q$, according to lemma
\ref{rep properties} there is a word $w$ representing it in the free
group $F_{\mathcal M}$ whose Magnus expansion $M(w)$ is an element
of the subring $S(C)$ of $R[C]$. Denote by $\overline{M}(w)$ the
image of $M(w)$ under the projection $S(C)\longrightarrow 1+Q(C)$,
so $\overline{M}(w)=p_2(p_1(M(w)))$ in the notation of
(\ref{expansion}).

\subsection{} {\em  Definition of ${\Phi}$ in the height $=1$ case.}
First consider the
special case when the first stage planar surface $P$ is a pair of
pants, $\partial P={\gamma}\cup{\alpha}_1\cup{\alpha}_2$. We will
follow the notation of section
\ref{presentation section}, and we will use the Magnus expansion
(\ref{Magnusexpansion}). In particular, the set $X$ of the variables
corresponding to the meridians to the handles of $C$ in $D^4$ is
divided into two subsets $X_{I_1}, X_{I_2}$, where the indices
reflect the components of the links $L_i\subset {\alpha}_i\times
D^2$ that the handles are attached to.

Let $Y_i$ be a monomial with non-repeating variables of maximal
length in the variables $X_{I_i}$, $i=1,2$, respecting the preferred
order (see \ref{preferred planar}). Note that $Q(C)$ in this case is
$2$-dimensional, spanned by the monomials $Y_1$, $Y_2$. Denoting by
$W_i$ a word representing the curve ${\wedge}_i$ in the free group,
given by the commutative diagram \ref{comm diagram}, note that
$\overline{M}(W_i)={\mu}_i\, Y_i$, where ${\mu}_i\neq 0$, $i=1,2$.

\begin{proposition} \label{pants}
\sl Given an element $g\in {\pi}_1(S^3\smallsetminus
L)/{\pi}_1(S^3\smallsetminus L)^q$, let $w$ be a word representing
it as in \ref{special word}, and consider its expansion in $1+Q(C)$:
$$\overline{M}(w)=1+{\alpha}_1 Y_1+{\alpha}_2 Y_2$$

for some ${\alpha}_1$, ${\alpha}_2$. Then ${\Phi}(g)\!
:={\mu}_2{\alpha}_1+{\mu}_1{\alpha}_2\in {\mathbb Z}$ is an invariant of $g$.
\end{proposition}

{\em Proof.} The coefficients ${\alpha}_i$ are well-defined with
respect to the relations $(R_1)-(R_3)$, see the discussion following
equation (\ref{expansion}). The relation $(R_4)$ is given by
$(W_1)^g(W_2)^{-1}$, and its expansion is of the form
$$\overline{M}((W_1)^g(W_2)^{-1})=1+{\mu}_1 Y_1-{\mu}_2 Y_2.$$

Let $w'$ be $w$ multiplied by a conjugate of $(W_1)^g(W_2)^{-1}$,
$\overline{M}(w')=1+{\alpha}'_1 Y_1+{\alpha}'_2 Y_2$. Then
${\alpha}'_1={\alpha}_1+{\mu}_1$, ${\alpha}'_2={\alpha}_2-{\mu}_2.$
Therefore ${\Phi}(w')={\Phi}(w)$. \qed

\smallskip

Consider the {\em general height $1$ case}: $\partial
P={\gamma}\cup{\alpha}_1\cup\ldots\cup{\alpha}_k$. As above, let
$Y_j$ be the preferred monomial in the variables $X_{I_j}$, and
$\overline{M}(W_j)=1+{\mu}_j Y_j$, ${\mu}_j\neq 0$, $j=1,\ldots,k$.
Define ${\mu}'_j=\prod_{i\neq j}{\mu}_i.$ The proof of the following
statement is a direct generalization of the proof in the pair of
pants case.

\smallskip

\begin{proposition} \label{height 1 invariant} \sl
Given an element $g\in {\pi}_1(S^3\smallsetminus
L)/{\pi}_1(S^3\smallsetminus L)^q$, as in proposition \ref{pants}
consider the expansion in $1+Q(C)$:
$\overline{M}(w)=1+\sum_j{\alpha}_j Y_j.$ Then ${\Phi}(g)\! :=\sum_j
{\alpha}_j{\mu}'_j$ is an invariant of $g$.
\end{proposition}

{\em Remark.} In fact there is a collection of $I_1 !\cdots I_k !$
invariants $\Phi$, parametrized by the monomials in non-repeating
variables $X_{I_1}, \ldots, X_{I_k}$. We chose a specific ${\Phi}$,
reflecting a particular choice of non-trivial $\bar\mu$-invariants of the homotopically
essential links $\widehat L_j$.

\subsection{} {\em Definition of the invariant ${\Phi}$ in the
general case.} \label{general phi} The definition is inductive.
Suppose the homomorphism ${\Phi}\co (Q(C),+)\longrightarrow
({\mathbb{Z}},+)$ is defined for $b$-cells of height $<h$, and let
$(C,{\gamma})$ be a $b$-cell of height $h$. $C$ is obtained from
$\overline{P}=P\times D^2$ by attaching $b$-cells $\{C_j\}_{j\in
I_i}$ of height $h-1$ to the components of links $L_i$, $L_i\subset
{\alpha}_i\times D^2$. Here $I_i$ is the (ordered) index set for the
components of $L_i$. As above, let ${\mu}_i$ be the non-trivial
$\bar\mu$-invariant of $\widehat L_i$ in the expansion of
${\wedge}_i$, with the given order of the components of $L_i$. Let
${\Phi}_j\co Q(C_j)\longrightarrow {\mathbb Z}$ denote the
inductively defined invariant of $C_j$. Recall from (\ref{tensor})
that
$$Q(C)=\oplus_i\;  \otimes_{j\in I_i} \; Q(C_j).$$

Denoting ${\mu}'_j=\prod_{i\neq j}{\mu}_i$, define $${\Phi}\co
Q(C)\longrightarrow {\mathbb Z} \; \; {\rm by} \; \;
{\Phi}=\sum_i{\mu}'_i (\otimes_{j\in I_i} {\Phi}_j).$$

\begin{proposition} \label{proposition phi} \sl
Given $g\in {\pi}_1(S^3\smallsetminus K)/({\pi}_1(S^3\smallsetminus
K))^q$, let $w$ be a word representing it in the free group, as in
\ref{special word}. Then ${\Phi}(\overline{M}(w))$ is well-defined,
and will be denoted ${\phi}(g)$.
\end{proposition}

{\em Proof.} The proof is inductive. The statement is true for
$b$-cells of height $1$ by proposition \ref{height 1 invariant}.
Suppose the statement is true for $b$-cells of height $<h$, and let
$C$ be a $b$-cell of height $h$. Assembling $C$ from $b$-cells of
height $h-1$ will be separated into two steps: (1) attaching them to
a link in a solid torus, and (2) attaching the results of step (1)
to a (planar surface)$\times D^2$, see section \ref{inductive} and
figures \ref{step 1}, \ref{step 2} in section \ref{representations}.

{\sl Step (1)}. Consider a $b$-cell $C$ of height $h$ such that
$C=S^1\times D^2\times I\cup (C_1\cup C_2)$ where the $b$-cells
$C_i$ have height $h-1$ and are attached along the components of a
link $L=(l_1,l_2)\subset S^1\times D^2\times\{ 1\}$. For simplicity
of notation, we assume $L$ has two components; the proof for a
larger number of components is directly analogous. Given a relation $r$ of type
$(R_4)$, let $I$ denote the ideal in $R[C]$ generated by the Magnus
expansion $M(W)-1$, where $W$ is a word representing $r$. It suffices to
prove that the intersection $I\cap Q(C)$ is in the kernel of
${\Phi}\co Q(C)\longrightarrow {\mathbb Z}$. The representation
$Q=Q(C)$ decomposes as $Q_1\otimes Q_2$ where $Q_i=Q(C_i)$, and
${\Phi}={\Phi}_1\otimes {\Phi}_2\co Q\longrightarrow {\mathbb Z}$,
so $$ker\, {\Phi}=(ker\, {\Phi}_1)\otimes Q_2+Q_1\otimes (ker\,
{\Phi}_2).$$

Since the bottom stage surface of $C$ is the annulus, there are
no relations $(R_4)$ at height $1$.  Therefore the relation $r$
corresponds to a body surface in either $C_1$ or $C_2$, say in
$C_1$.

First we impose an {\em additional assumption} that, in the context
of definition \ref{admissible links}, for each link $L$ defining the
$b$-cell $C$ there is a word $W$ representing $\wedge$ in the free
group such that $W$ involves only the variables $m_1,\ldots,m_n$,
and not the longitude $l$ of the solid torus. For example, this
assumption is satisfied in the central case $L=$(iterated) Bing double.
After giving a proof in this restricted setting, we show how the
argument goes through in the general case. The assumption above
implies that each relation $r$ of type $(R_4)$ has a word
representing it in the free group, whose Magnus expansion is an
element of either $R[C_1]$ or $R[C_2]$.

Let $r\in R[C_1]\subset R[C]$ be a relation, and denote by $I_1$ and
$I$ the ideals generated by $r$ in $R[C_1]$, $R[C]$ respectively.
Observe that $I\cap Q(C)=I\cap(Q(C_1)\otimes Q(C_2))=(I\cap
Q(C_1))\otimes Q(C_2)$. Since $I_1\subset ker\, {\Phi}_1$, $I\cap
Q(C)\subset ker\, {\Phi}$, and the proof is complete.

Now consider the general case, i.e. we remove the extra assumption
imposed in the paragraph above. The difference with that case is
that even though $r$ is a relation corresponding to $C_1$, one
cannot assume that $r$ is an element of the subring $R[C_1]$ of
$R[C]$. However (see end of section \ref{links in torus}) $\wedge$
has a word representing it whose expansion is of the form
$1+x_{i_1}\cdots x_{i_n}+$higher order terms. That is, all {\em
first non-vanishing terms} with non-repeating variables in its
Magnus expansion are elements of $R[C_1]$. The proof is completed by
the observation that only first non-vanishing terms contribute to
$I\cap Q(C)$.

{\sl Step (2)}, see figure \ref{step 2}. Now $C$ equals $(P\times
D^2)\cup C_1\cup C_2$, where $P$ is a planar surface, the $b$-cells $C_i$ have height $h$ and
whose bottom stage surfaces are annuli. For simplicity of notation
we assume $P$ is a pair of pants; the case of a planar surface with
more boundary components is treated analogously. Denoting $\partial
P={\gamma}\cup{\alpha}_1\cup{\alpha}_2$, $C_i$ is attached along
${\alpha}_i\times D_2$, $i=1,2$. In this case $$R[C]\cong
R[C_1]\oplus R[C_2],\; \, Q(C)\cong Q(C_1)\oplus Q(C_2).$$ As above,
given a relation $r$ of type $(R_4)$, we need to show $$I\cap
Q(C)\subset ker({\Phi}\co Q(C)\longrightarrow {\mathbb Z}).$$ We
have ${\Phi}={\mu}_2{\Phi}_1\oplus {\mu_1}{\Phi}_2$. There are two
cases to consider: $r$ corresponds to a surface in $C$ at height
$>1$, or it is a new relation corresponding to $P$. In the first
case, one may assume $r\in R[C_1]$. Denote by $I_1$, $I$ the ideals
generated by $r$ in $R[C_1]$, $R[C]$. Then since $R[C]$ is a direct
sum of rings $R[C_1]\oplus R[C_2]$, $I=I_1\subset R[C_1]\subset
R[C]$. Clearly then $I\subset ker({\Phi})$.

Consider the second case: $r$ is a new relation, corresponding to
the bottom stage surface $P$ of $C$. Denote the meridians to $L_1$ by
$m'_1,\ldots, m'_k$ and the meridians to $L_2$ by
$m''_1,\ldots,m''_l$; let $\{ x'_i\}, \{ x''_j\}$ be the
corresponding variables. Then the Magnus expansion of $r$ is of the form
$$M(r)=1+{\mu}_1x'_1\cdots x'_k-{\mu}_2 x''_1\cdots x''_l+{\rm
higher}\;\; {\rm order}\;\;{\rm terms}. $$ Consider the image of $r$
in $R[C]$. Note that the first term ${\mu}_1x'_1\cdots x'_k$ is in
$R[C_1]$, the second term ${\mu}_2 x''_1\cdots x''_l$ is in
$R[C_2]$, and in fact all higher order terms vanish in $R[C]$, since
the first  non-vanishing terms already have maximal length. Any
element of $R[C]$ of the form ${\mu}_1Y+{\mu}_2 Z$, where $Y\in
R[C_1]$, $Z\in R[C_2]$, is in the kernel of ${\Phi}$. Therefore
$r\in ker({\Phi})$, and any other element in the ideal generated by
$r$ is longer and vanishes in $R[C]$ (so in fact $I=\{ r\} \subset
ker({\Phi})$.) This concludes the proof of proposition
\ref{proposition phi}. \qed

Proposition \ref{proposition phi} constructs a homomorphism
${\phi}\co {\pi}_1(S^3\smallsetminus K)\longrightarrow {\mathbb Z}$.
In particular, ${\phi}(g)$ is well-defined
with respect to multiplication by elements of the relation
subgroup in $F_{\mathcal M}$, so ${\phi}(g)\neq 0$ implies $g\neq
1\in {\pi}_1(D^4\smallsetminus C)/({\pi}_1(D^4\smallsetminus C))^q$.
It suffices to prove lemma \ref{invariant} for $g$ equal to a meridian
$m$ to the knot $K$ in $S^3$. The fact that ${\phi}(m)\neq 0$ is proved by inspection: at each
surface stage $P$ of $C$, $\partial P={\gamma}_P\cup_i{\alpha}_i$,
the meridian to $P$ is conjugate to the ${\wedge}$-curve
corresponding to the solid torus ${\alpha}_i\times D^2$, for any
given $i$. Applying the analysis at the end of section \ref{links in
torus} inductively to the meridians to the surface stages of $C$,
moving up from the meridian $m$ to the bottom stage, one observes
that there is a word $w$ representing $m$ in $F_{\mathcal M}$ such
that $\overline M(w)$ is a generating monomial for $Q(C)$. Due to
the tensor decompositions of $Q(C)$ and $\Phi$,
${\phi}(m)={\Phi}(\overline M(w))\neq 0$. This concludes the proof of
lemma \ref{invariant}. \qed

\section{Applications to link homotopy: proof of theorem \ref{homotopically trivial}} \label{homotopy}

This section shows how the theory of Bing cells fits in the
framework of Milnor's theory of link homotopy. We generalize the
invariant $\Phi$ defined in the previous section to a collection of
Bing cells to prove theorem \ref{homotopically trivial}.

{\em Proof of theorem} \ref{homotopically trivial}. Let
$L=(l_1,\ldots,l_n)$ and suppose the components of $L$ bound
disjoint Bing cells $C_1\ldots, C_n$ in $D^4$. Denote $C=\cup_i
C_i$. Suppose $L$ is homotopically essential, and without loss of
generality one may assume $L$ is almost homotopically trivial, so
there is a well-defined and non-trivial ${\mu}$-invariant with
non-repeating coefficients of length $n$. Order the components of
$L$ so that ${\mu}_{1\ldots n}(L)\neq 0$.

The results of sections \ref{invariants}, \ref{generalized Milnor},
\ref{representations} and \ref{invariant phi} generalize from the
setting of a single $b$-cell as follows. Let $\mathcal{M}_i$ denote
a set of meridians to the handles of $C_i$. By Alexander duality
$H_1(D^4\smallsetminus C)$ is generated by
$\mathcal{M}=\sqcup_i\mathcal{M}_i$. Denote the corresponding
variables for the Magnus expansion by $X_i$, $X=\sqcup X_i$. Again by
Alexander duality, the relations in ${\pi}_1(D^4\smallsetminus
C)/({\pi}_1(D^4\smallsetminus C))^q$ are all of types $(R_1)-(R_4)$
(see section \ref{presentation section}), contributed by the
$b$-cells $C_i$. Each relation of type $(R_1)-(R_3)$ involves only
variables in a single set $\mathcal{M}_i$. The assumptions on the links
defining the $b$-cells in section \ref{links in torus} imply that all first non-vanishing terms in
the Magnus expansion of any relation of type $(R_4)$ also involve
the variables in a single $X_i$. Variables from other sets $X_j$ may
be present, but only in higher-order terms.

Define $GM(C)$ as the free group $F_{\mathcal{M}}$ modulo relations
(\ref{Milnor relations}), where all of the meridians $m,m_1,m_2$
involved in the commutators in (\ref{Milnor relations}) are elements
of the same $\mathcal{M}_i$, for any given $1\leq i\leq n$. Define
$R[C]$ as the quotient of ${\mathbb Z}\{ X\}$ by the ideal
introduced in definition \ref{non-repeated} where the variables
$x_I$, $x_{I'}$ are elements of $X_i$ for the same $i$. Consider the
Magnus expansion in the following diagram, analogous to that in proposition
\ref{proposition diagram}:

\[
\xymatrix{
    F_{\mathcal{M}}
    \ar[r] \ar[d]
     &
    GM(F_{\mathcal{M}})
     \ar[r] \ar[d] &
    GM(C)
    \\
    {\mathbb Z}\{X\}
    \ar[r] &
    R[C]&
    }\]

Following definitions \ref{rep R}, \ref{rep S}, \ref{rep Q},
introduce $\widetilde R(C)$, $S(C)=1+\widetilde R(C)$. Define $Q(C)$
using the order on the components of $L$ reflecting a non-trivial
${\mu}$-invariant (see above):
$$Q(C)=Q(C_1)\otimes\ldots \otimes Q(C_n).$$

The proof of lemma \ref{rep properties} goes through, in particular
given any element $g\in {\pi}_1(S^3\smallsetminus
L)/({\pi}_1(S^3\smallsetminus L))^q$, there is a word $w_0$
representing it in $F_{\mathcal M}$ such that $M(w_0)\in S(C)$.
Denote by $\overline M$ the composition of the Magnus expansion $M$
with the projection $S(C)\longrightarrow 1+Q(C)$. Denoting by
${\Phi}_i$ the homomorphism $Q(C_i)\longrightarrow{\mathbb Z}$
defined in \ref{general phi}, consider
    $${\Phi}=\otimes_i {\Phi}_i\co Q(C)=\otimes_i Q(C_i)\longrightarrow {\mathbb
Z}.$$

    Given $g\in {\pi}_1(S^3\smallsetminus
L)/({\pi}_1(S^3\smallsetminus L))^q$, ${\Phi}(\overline M(w_0))$ is
a well-defined integer. Moreover, if ${\Phi}(\overline M(g))\neq 0$,
then $i_*(g)\neq 1\in \pi_1(D^4\smallsetminus
C)/(\pi_1(D^4\smallsetminus C))^q$. Consider the commutative diagram

\[
\begin{CD}
\pi_1(S^3\smallsetminus L)/(\pi_1(S^3\smallsetminus L))^q
@>{i_*}>>\pi_1(D^4\smallsetminus C)/(\pi_1(D^4\smallsetminus C))^q \\
@A{p_1}AA
@A{p_2}AA\\
F_{m_1,\ldots,m_{n}}@>{\alpha}>>
F_{\mathcal M}=F_{\mathcal{M}_1,\ldots,\mathcal{M}_n}\\
@V{M_1}VV @V{M}VV\\
{\mathbb{Z}}\{x_1,\ldots,x_n\} @>{\beta}>>
{\mathbb Z}\{X\}={\mathbb{Z}}\{X_1,\ldots,X_n\} \\
\end{CD}
\]

\smallskip

Recall from the proof of lemma \ref{invariant} at the end of section
\ref{invariant phi} that each meridian $m_i$ has a word $w_i$
representing it in $F_{\mathcal M}$ such that $\overline M(w_i)$ is
a generating monomial for $Q(C_i)$, and ${\Phi}_i(\overline
M(w_i))\neq 0$. In the diagram above ${\alpha}$ is defined by
setting ${\alpha}(m_i)=w_i$. Then ${\beta}$ is given by
${\beta}(x_j)=M({\alpha}(m_j))-1$.

Since $L$ is homotopically essential, there is a relation
$[m_i,l_i]$ in $\pi_1(S^3\smallsetminus L)/(\pi_1(S^3\smallsetminus
L))^q$ such that the Magnus expansion $M_1$ of a word $W$
representing it in $F_{m_1,\ldots,m_n}$ is of the form
$1+{\mu}x_1\cdots x_n+\ldots$ where ${\mu}\neq 0$. However the
projection of ${\beta}(x_1\cdots x_n)$ onto $Q(C)$ is a product of
generating monomials, one for each $Q(C_i)$, and it follows from the
definition of $\Phi$ that ${\Phi}({\alpha}(W))\neq 0$. Since
${\Phi}(\overline M(w_0))$ is an invariant of $g\in
\pi_1(S^3\smallsetminus L)/(\pi_1(S^3\smallsetminus L))^q$, where
$p_1(w_0)=g$, $p_1(W)\neq 1\in \pi_1(S^3\smallsetminus
L)/(\pi_1(S^3\smallsetminus L))^q$. But $p_1(W)=[m_i,l_i]$ is a
relation in that group. This contradiction concludes the proof of
theorem \ref{homotopically trivial}. \qed


\bigskip


\begin{thebibliography}{10}

\bibitem{Cochran} T. Cochran,
{\em Derivatives of links: Milnor's concordance invariants and Massey's products},
Mem. Amer. Math. Soc. 84 (1990), no. 427.

\bibitem{Dwyer} W. Dwyer, {\em Homology, Massey products and maps
between groups}, J. Pure Appl. Algebra 6 (1975), 177-190.

\bibitem{Freedman2} M.H. Freedman, {\em Are the Borromean rings
$(A,B)$-slice?}, Topology Appl., 24 (1986), 143-145.


\bibitem{FK} M. Freedman and V. Krushkal, {\em Topological arbiters}, arXiv:1002.1063, to appear in J. Topol.

\bibitem{FL} M. Freedman and X.S. Lin, {\em On the $(A,B)$-slice
problem}, Topology Vol. 28 (1989), 91-110.

\bibitem{FQ} M. Freedman and F. Quinn, {\em The topology of
4-manifolds}, Princeton Math. Series 39, Princeton, NJ, 1990.

\bibitem{Giffen} C. Giffen, {\em Link concordance implies link homotopy}, Math. Scand.
45 (1979), 243-254.


\bibitem{Goldsmith} D. Goldsmith,
{\em Concordance implies homotopy for classical links in $M^3$},
Comment. Math. Helv. 54 (1979), 347-355.

\bibitem{K} V. Krushkal, {\em Additivity properties of Milnor's
$\bar\mu$-invariants}, J. Knot Theory Ramifications 7 (1998),
625-637.


\bibitem{K2} V. Krushkal, {\em Link groups and the A-B slice problem},
Topology and physics, 220-235, Nankai Tracts Math.
12, World Sci. Publ., Hackensack, NJ, 2008 [arXiv:math/0602105]

\bibitem{K3} V. Krushkal, {\em A counterexample to the strong version of 
Freedman's conjecture}, Ann. of Math. 168 (2008), 675-693 [arXiv:math/0610865]

\bibitem{KT} V. Krushkal and P. Teichner, {\em Alexander duality,
gropes and link homotopy}, Geom. Topol. 1 (1997), 51-69 [arXiv:math/9705222]


\bibitem{Milnor1} J. Milnor, {\em Link Groups}, Ann. Math 59 (1954), 177-195.

\bibitem{Milnor2} J. Milnor, {\em Isotopy of links}, Algebraic geometry and
topology, Princeton Univ. Press, 1957, pp. 280-306.

\bibitem{Stallings} J. Stallings, {\em Homology and central series of groups},
J. Algebra 2 (1965), 1970-1981.

\bibitem{Warfield} R. Warfield, {\em Nilpotent groups},
Lecture Notes in Mathematics, Vol. 513,
Springer-Verlag, Berlin-New York, 1976.


\end{thebibliography}
\end{document}